\documentstyle[12pt]{amsart}

\newcommand{\catqot}{/\hskip-3pt/}
\newcommand{\C}{{\Bbb C}}

\newcommand{\E}{{\cal E}}

\newcommand{\F}{{\cal F}}
\newcommand{\GL}{\mathop{\rm GL}}

\newcommand{\Hom}{\mathop{\rm Hom}}

\renewcommand{\Im}{\mathop{\rm Im}}
\newcommand{\id}{\mathop{\rm id}}

\renewcommand{\L}{{\cal L}}

\newcommand{\n}{{\cal N}}
\renewcommand{\O}{{\cal O}}
\renewcommand{\P}{{\Bbb P}}
\newcommand{\Pic}{\mathop{\rm Pic}}
\newcommand{\Proj}{\mathop{\rm Proj}}
\newcommand{\Q}{{\Bbb Q}}

\newcommand{\SL}{\mathop{\rm SL}}
\newcommand{\Spec}{\mathop{\rm Spec}}
\renewcommand{\tilde}{\widetilde}
\newcommand{\Z}{{\Bbb Z}}

\newcommand{\la}{\lambda}
\newcommand{\lra}{\longrightarrow}
\newcommand{\ra}{\rightarrow}

\newcommand{\lma}{\longmapsto}
\newcommand{\p}{\prime}
\newcommand{\q}{\quad}
\newcommand{\Quot}{\mathop{\rm Quot}}

\renewcommand{\phi}{\varphi}
\newcommand{\rk}{\mathop{\rm rk}}
\newcommand{\eps}{\varepsilon}
\theoremstyle{plain}
\newtheorem{Thm}{Theorem}[subsection]
\newtheorem{Cor}[Thm]{Corollary}
\newtheorem{Prop}[Thm]{Proposition}
\newtheorem{Lem}[Thm]{Lemma}

\newtheorem{Claim}{\rm Claim}

\theoremstyle{remark}
\newtheorem{Rem}[Thm]{\bf Remark}
\newtheorem{Ex}[Thm]{\bf Example}

\title{Master spaces for stable pairs}
\author{Ch.\ Okonek, A.\ Schmitt, A.\ Teleman}
\address
{\hskip -4mm {Institut f\"ur Mathematik\\ Universit\"at Z\"urich\\
Winterthurerstr.\ 190\\ \hfill{\break} CH-8057 Z\"urich,  Switzerland}}
\thanks{The authors were partially supported by AGE --- Algebraic Geometry
in Europe,
contract No.\ ER-BCHRXCT 940557 (BBW 93.0187), and by SNF, Nr.\ 2000 --
045209.95/1}

\begin{document}
\maketitle
\centerline{$\underline{\hskip 5truecm}$}
\vskip 1.2truecm

\def\map{\longrightarrow}
\def\textmap#1{\mathop{\vbox{\ialign{
                                ##\crcr
    ${\scriptstyle\hfil\;\;#1\;\;\hfil}$\crcr
    \noalign{\kern-1pt\nointerlineskip}
    \rightarrowfill\crcr}}\;}}

\section*{Introduction}

In this paper we   construct master spaces for certain coupled vector
bundle problems over a fixed projective variety $X$.

From a technical point of view, master spaces classify oriented pairs $({\cal
E},\varepsilon,\varphi)$ consisting of a torsion  free coherent sheaf
${\cal E}$ with
fixed Hilbert polynomial, an orientation $\varepsilon$ of the determinant
of ${\cal E}$,
and a framing $\varphi:{\cal E}\map {\cal E}_0$ with values in a fixed
reference sheaf
${\cal E}_0$, satisfying certain semistability conditions. The relevant
stability concept
is new and does not involve the choice of a parameter, but it can easily be
compared to
the older parameter-dependent  stability concepts for (unoriented)  pairs.

The corresponding moduli spaces ${\cal M}$ have the structure of polarized
projective
varieties endowed with a natural $\C^*$-action which can be exploited in
two interesting
ways:

1. The fixed point set ${\cal M}^{\C^*}$ of the $\C^*$-action is a union
$${\cal
M}^{\C^*}={\cal M}_{source}\cup{\cal M}_{sink}\cup {\cal M}_R\ ,$$
 where ${\cal
M}_{source}$ is a Gieseker moduli space of semistable oriented sheaves,
${\cal M}_{sink}$
is a certain (possible empty) Grothendieck Quot-scheme, and the third term
${\cal
M}_R:={\cal M}^{\C^*}\setminus({\cal M}_{source}\cup{\cal M}_{sink})$ is
the so-called
"variety of reductions", which consists essentially of lower rank objects. The structure as
a $\C^*$-space can be used to relate "correlation functions" associated
with the different
parts of ${\cal M}^{\C^*}$ to each other [OT2].

2. Master spaces are also useful for the investigation of the birational
geometry of
the moduli spaces ${\cal M}_\delta$ of $\delta$-semistable pairs in the
sense of [HL2].

Indeed, each of the ${\cal M}_\delta$'s can be obtained as a suitable
$\C^*$-quotient of
the master space ${\cal M}$, and it can be shown that every two quotients
${\cal M}_\delta$,
${\cal M}_{\delta'}$ are related by a  chain of generalized flips in the
sense of [Th].

When $X$ is a projective curve with trivial reference sheaf ${\cal E}_0={\cal
O}_X^{\oplus k}$, our master space can be considered as a natural
compactification of the
one described in [BDW]. Their space becomes an open subset of ours whose
complement is
the Quot-scheme ${\cal M}_{sink}$ alluded to above (${\cal M}_{sink}$ is
empty iff
$k<\rk({\cal E})$). Applying the ideas of 1.\ in this situation leads to
formulas for
volumina and characteristic numbers and to a new proof of the Verlinde
formula when $k=1$,
and allows to relate Gromov-Witten invariants for Grassmannians to simpler
vector bundle
data when
$k>\rk({\cal E})$.

In the case of an algebraic surface $X$, master spaces can be viewed as
algebraic analoga
of certain gauge theoretic  moduli spaces of monopoles which can be used to
relate
Seiberg-Witten invariants and Donaldson polynomials [OT1], [T1]. The latter
application
was actually our original motivation for the  construction of master spaces.

The study of non-abelian monopoles on K\"ahler surfaces leads naturally to the
investigation of a certain moment map on an infinite dimensional K\"ahler space.

The associated stability concept, which is expected to exist on general
grounds [MFK], is
precisely the one which gave rise to the stability definition for oriented
pairs [OT2].
Since the moduli space of non-abelian monopoles admit an Uhlenbeck type
compactification
[T1], it was natural to look for a corresponding Gieseker type
compactification of their
algebro-geometric analoga. These compactifications, the master spaces for
stable pairs,
provide   very useful models for understanding the ends of monopole moduli
spaces in the more difficult gauge theoretical context [T2]. Understanding
these ends is the
essential final step in our program for relating Donaldson polynomials and
Seiberg-Witten
invariants [OT1], [T1]. Let us now briefly describe the main ideas and
results of this
paper.

 The construction
of master spaces requires the study of GIT-quotients for direct sums of
representations, i.e.\ the construction of quotients  $\P(A\oplus
B)^{ss}\catqot G$, where $G$ is a reductive group acting linearly on vector
spaces $A$ and $B$. Since the Hilbert criterion is difficult to apply in
this situation, we have chosen another approach instead. The idea is to use
the natural $\C^*$-action $z\cdot\langle a,b\rangle:=\langle a,zb\rangle$ on
$\P(A\oplus B)$ which commutes with the given action of $G$. Our first main
result characterizes $G$-semistable points in $\P(A\oplus B)$ in terms of
$G$-semistability of their images in all possible
$\C^*$-quotients of $\P(A\oplus B)$. The proof is based on a commuting
principle for actions of products of groups.

These results, which we prove in the first section, explain in particular
why chains of flips occur in GIT-problems for $G\times\C^*$-actions [DH], [Th].

In the second section of our paper, after defining stability for oriented pairs
$({\cal E},\varepsilon,\varphi)$, we prove a crucial boundedness result  and
construct the corresponding parameter space ${\frak B}$. This space admits
a morphism
$\iota:{\frak B}\map \P({\frak Z})$ into a certain Gieseker space
$\P({\frak Z})$ which is equivariant w.r.t.\ a natural action of a product
$\SL\times\C^*$ of $\C^*$ with a special linear group. The $\SL$-action  on
$\P({\frak Z})$ possesses a linearization in a suitable line bundle, and the
preimage of the subset $\P({\frak Z})^{ss}$ of $\SL$-semistable points is
precisely the open subspace ${\frak B}^{ss}\subset {\frak B}$ of points
representing semistable oriented pairs. In order to prove this, we apply
our GIT-Theorem from the first section to the $\SL\times\C^*$-action on
$\P({\frak Z})$, and thereby reduce the proof to results in [G] and [HL1].

Then we show that the induced map $\iota|_{{\frak B}^{ss}}:{\frak B}^{ss}\map
\P({\frak Z})^{ss}$ is finite and hence descends to a finite map
$\bar\iota:{\frak B}^{ss}\catqot \SL\map \P({\frak Z})^{ss}\catqot \SL$. The
quotient  ${\frak B}^{ss}\catqot \SL$, which is therefore a projective
variety, is our master space.

The  ideas and techniques of this paper can also be applied
to construct master spaces in  other interesting  situations, e.g.\ by
coupling with
sections in twisted endomorphism bundles. When $X$ is a curve and the twisting
line bundle
is the canonical bundle, one obtains a natural compactification of the moduli
 spaces of
Higgs bundles [H], [S].

Similar ideas should also apply to coupling with singular
objects like parabolic structures. We refer to [OT2] for a general
description of
the  underlying coupling principle and its application to computations of
correlation
functions.

%\end{document}
%%%%%%%%%%%%%%%%%%%%%%%%%%%%%%%%%%%%%%%%%%%%%%%%%%%%%%%%%%%

\subsection*{Conventions}

Our ground field is $\C$. A \it polarization \rm on a quasi-projective
variety $X$ is
an equivalence class $[L]$ of ample line bundles,
where two line bundles $L_1$ and $L_2$ are \it equivalent\rm , if there
exist positive integers $n_1$ and $n_2$ such that
$L_1^{\otimes n_1}\cong {L_2}^{\otimes n_2}$.
\par
If $W$ is a finite dimensional vector space, we denote by $\P(W)$ its
projectivization in the sense of Grothendieck, i.e., the  closed
points of $\P(W)$
correspond to
lines in the dual space $W^\vee$.
We do not distinguish notationally between a vector space $W$ and its
 associated scheme.

\section{A theorem from Geometric Invariant Theory}
\subsection{Background material from GIT}
\label{BackGIT}
Let $G$ be a reductive algebraic group and let
$\gamma\colon G\lra \GL(W)$
be a rational representation in the finite dimensional vector space $W$.
The map $\gamma$ defines an action of $G$ on the dual space $W^\vee$
given by
$$g\cdot w:= w\circ \gamma(g^{-1})\qquad \forall g\in G; w\in W^\vee,$$
an action $\overline{\gamma}$ on the projective space $\P(W)$,
and a linearization of this action in $\O_{\P(W)}(1)$.
In the following we identify $H^0(\P(W),\O_{\P(W)}(k))$ with
$S^kW$.\par
Recall that a point $x\in \P(W)$ is \it $\gamma$-semistable \rm
if and only if the orbit closure $\overline{G\cdot w}$
of any lift $w\in W^\vee\backslash\{0\}$ does not contain $0$.
Denote by $\P(W)_{\gamma}^{ss}\subset\P(W)$ the open set of semistable
points and by $\P(W)_\gamma^{ps}$ the set of \it $\gamma$-polystable \rm
points, i.e.\ the semistable points whose orbit is closed in
$\P(W)_{\gamma}^{ss}$.
Equivalently, a point $x\in \P(W)$ is polystable if and only if the
orbit $G\cdot w$ of any lift $w\in W^\vee\backslash\{0\}$
is closed in $W^\vee$.
With this terminology, $x\in \P(W)$ is \it $\gamma$-stable \rm if and only if it
is polystable and its stabilizer $G_x$ is finite.
Let $\pi_\gamma\colon \P(W)_\gamma^{ss}\lra Q_\gamma:= \P(W)\catqot_\gamma G$
be the categorical quotient.
For sufficiently large $n$, $Q_\gamma$ admits a projective
embedding
$j_n\colon Q_\gamma\hookrightarrow \P({S^nW}^G)$
such that the following diagram commutes:
%%%%%%%%%%%%%%%%%%%%%%%%%%%%%%%%%%%
\begin{equation}
\label{eqemb}
\begin{array}{c}
\unitlength=1mm
\begin{picture}(70,24)(0,9)
\put(0,29){$\P(W)^{ss}_\gamma\subset \P(W)$}
\put(34,30){\vector(1,0){20}}
\put(34,31){\oval(3,1.8)[l]}
\put(40,32){${\scriptstyle v_n}$}
\put(56,29){$\P(S^n W)\ $}
\put(7,25){\vector(0,-1){10}}
\multiput(64,24)(0,-2){4}{\line(0,1){1}}
\put(64,16){\vector(0,-1){1}}
\put(5,9){$Q_\gamma$}
\put(14,11){\oval(3,1.8)[l]}
\put(14,10){\vector(1,0){40}}
\put(56,9){$\P(S^nW^G)$}
\put(2,20){$\scriptstyle\pi_\gamma$}
\put(38,12){$\scriptstyle j_n$}
\put(66,20){${\scriptstyle p_G}$}
\end{picture}
\end{array}
\end{equation}
%%%%%%%%%%%%%%%%%%%%%%%%%%%%%%%%

In this diagram, $v_n$ stands for the $n$-th Veronese embedding and
$p_G$ is the projection induced by the inclusion ${S^nW}^G\subset
S^nW$.
The space $Q_\gamma$ comes with a natural polarization represented by
$L_n:= j_n^*\O_{\P({S^nW}^G)}(1)$.
Indeed, by (\ref{eqemb})
we have $\pi_\gamma^*L_n\cong\O_{\P(W)_\gamma^{ss}}(n)$,
and from the commutative diagram
%%%%%%%%%%%%%%%%%%%%%%%%%%%%%%%%%%%%%%%%
\begin{equation}
\label{eqemb}
\begin{array}{c}
\unitlength=1mm
\begin{picture}(120,24)(0,9)
\put(0,29){$\P(W)^{ss}_\gamma\subset \P(W)$}
\put(34,30){\vector(1,0){10}}
\put(34,31){\oval(3,1.8)[l]}
\put(36,32){${\scriptstyle v_{n_1}}$}
\put(46,29){$\P(S^{n_1} W)$}
\put(7,25){\vector(0,-1){10}}
\put(5,9){$Q_\gamma$}
\put(14,11){\oval(3,1.8)[l]}
\put(14,10){\vector(1,0){15}}
\put(31,9){$\P(S^{n_1}W^G)$}
\put(2,20){$\scriptstyle\pi_\gamma$}
\put(18,12){$\scriptstyle j_{n_1}$}
\put(110,20){${\scriptstyle p_G}$}
\put(53,10){\vector(1,0){10}}
\put(53,11){\oval(3,1.8)[l]}
\put(55,12){$\scriptstyle v_{n_2}$}
\put(66,9){$\P(S^{n_2}(S^{n_1}W^G))$}
\multiput(116,24)(0,-2){4}{\line(0,1){1}}
\put(116,16){\vector(0,-1){1}}
\multiput(97,10)(2,0){3}{\line(1,0){1}}
\put(104,10){\vector(1,0){1}}
\put(107,9){$\P(S^{n_1n_2} W^G)$}
\put(67,30){\vector(1,0){36}}
\put(67,31){\oval(3,1.8)[l]}
\put(80,32){${\scriptstyle v_{n_2}}$}
\put(106,29){$\P(S^{n_1n_2} W)$}
\end{picture}
\end{array}
\end{equation}
%%%%%%%%%%%%%%%%%%%%%%%%%%%%%%%%%%%%%%%%%%%%%%%%%%%

we infer $L_{n_1}^{\otimes n_2}\cong L_{n_1n_2}$, hence
\begin{equation}
\label{polind}
L_{n_1}^{\otimes n_2}\cong L_{n_2}^{\otimes n_1},\q \forall n_1,n_2\q
\hbox{large enough}.
\end{equation}
\begin{Rem}
\label{MoreGeneral}
In the following, we will mainly consider actions on projective spaces.
However, if $X$ is a quasi-projective variety with an action of an
algebraic group $G$ which is
linearized in an ample line bundle $L$, then $L^{\otimes n}$ induces,
for $n$ large enough, a $G$-invariant embedding of $X$
into $\P:=\P(H^0(L^{\otimes n}))$ such that the semistable,
polystable, and stable points of $X$ are mapped to the semistable, polystable,
and stable points of $\P$.
Hence all the results which we will prove hold also in this more general
setting,
and will be used in this generality in Section 2.
\end{Rem}

\subsection{Polarized $\C^*$-quotients}
\label{PolC*Quot}
Let $\la\colon\C^*\ra \GL(W)$ be a rational representation of $\C^*$ in the
finite dimensional vector space $W$ and let
$\overline{\la}\colon \C^*\times\P(W)\lra \P(W)$ be the induced action.
The space $W^\vee$ splits as a direct sum
$$W^\vee=\bigoplus_{i=1}^m W^\vee_i,$$
where $W^\vee_i$ is the eigenspace of the character
$\chi_{d_i}\colon\C^*\lra\C^*, z\lma z^{d_i}$.
We assume $d_1<d_2<\cdots<d_m$.
Let $x\in \P(W)$ and choose a lift $w\in W^\vee\backslash\{0\}$ of $x$.
Define
\begin{eqnarray*}
d^\la_{\min}(x)&:=&\min\bigl\{\, d_i\ \vert\ w\hbox{ has a non-trivial
component in } W_i^\vee\,\bigr\}\\
d^\la_{\max}(x)&:=&\max\bigl\{\, d_i\ \vert\ w\hbox{ has a non-trivial
component in } W_i^\vee\,\bigr\} .
\end{eqnarray*}

\begin{Prop}
\label{C^*-ss}\hfill{\break}
{\rm i)} A point $x\in\P(W)$ is $\la$-semistable if and only if
$d^\la_{\min}(x)\le 0\le d^\la_{\max}(x)$.\\
{\rm ii)} A point $x\in\P(W)$ is $\la$-polystable if and only if
either $d^\la_{\min}(x)=0=d^\la_{\max}(x)$ or
$d^\la_{\min}(x)<0<d^\la_{\max}(x)$.
\end{Prop}

\begin{pf}
Let $w=(w_1,...,w_n)\in W^\vee\backslash\{0\}$ be a lift of $x$, where we take
coordinates with respect to a basis of eigenvectors.
For $z\in \C^*$, we get
$$z\cdot w=(0,...,0,z^{d^\la_{\min}(x)}\cdot
w_{i_0},...,z^{d^\la_{\max}(x)}\cdot
w_{i_r},0,...,0).$$ Using this description, the assertion becomes obvious.
\end{pf}

As remarked above, we can view $\la$ as a linearization
of the action $\overline{\la}$.
There are two natural ways of changing this linearization:
\begin{enumerate}
\item Multiplying $\la$ by a character: Let $d$ be an integer, and denote by
$\la_d$ the representation
$z\lma z^d\cdot \la(z)$ of $\C^*$ in $\GL(W)$. This means that we change the
$\O_{\P(W)}(1)$-linearization
of $\overline{\la}$ by multiplying it with the character $\chi_{d}\colon
\C^*\lra\C^*, z\lma z^{d}$.
\item Replacing $\la$ by a symmetric power:
Let $\la^k\colon\C^*\lra \GL(S^kW)$ be the $k$-th symmetric power of $\la$.
This induces an $\O_{\P(W)}(k)$-linearization of $\overline{\la}$.
\end{enumerate}
Now we combine both methods, i.e.,
we change $\la^k$ to the representation $\la^k_d$ of $\C^*$
in $\GL(S^kW)$. As above, this defines an $\O_{\P(W)}(k)$-linearization of
$\overline{\la}$.
Altogether, we have a family $\la_d^k$, $k\in\Z_{>0}$, $d\in\Z$,
of linearizations of $\overline{\la}$.
Since two $\O_{\P(W)}(k)$-linearizations of $\overline{\la}$ differ
by a character of $\C^*$, these are indeed all possible
linearizations.\par
Every linearization $\la_d^k$ yields a polarized
GIT-quotient $\bigl(Q_d^k:=\P(W)\catqot_{\la_d^k}\C^*, [L^k_d]\bigr)$, and
$(Q_d^k,[L^k_d])$ and $(Q_{d^\p}^{k^\p},[L^{k^\p}_{d^\p}])$ are isomorphic
as polarized varieties
when the ratios $d/ k$ and $d^\p/ k^\p$ coincide.
To see this, one just has to observe that, for any positve integer $t$, the
linearization $\la_{t\cdot d}^{t\cdot k}$ is the $t$-th
symmetric power of the linearization $\la_d^k$.
\par
Since for a point $x\in \P(W)$ we have
$$d^{\la^k_d}_{\min}(x)=k\cdot d^\la_{\min}-d,\q
d^{\la_d^k}_{\max}(x)=k\cdot d^{\la}_{\max}-d,$$
we obtain the following corollary to Proposition~\ref{C^*-ss}:

\begin{Prop}
\label{C^*-ss2} \hfill{\break}
{\rm i)} The point $x$ is $\la_d^k$-semistable if and only if
$d^{\la}_{\min}(x)\le d/ k\le d^\la_{\max}(x)$.
\\
{\rm ii)} The point $x$ is
$\la_d^k$-polystable if and only if either
$d^{\la}_{\min}(x)=d/ k= d^\la_{\max}(x)$ or
$d^{\la}_{\min}(x)< d/ k< d^\la_{\max}(x)$.
In particular, every point $x\in\P(W)$ is $\la_d^k$-polystable
for suitable numbers
$k\in\Z_{>0}$, $d\in\Z$.
\end{Prop}

For integers $i$ with $1\le i\le 2m$ we define the following intervals in
$\P^1_\Q$:
$$I_i:=\cases  \P_\Q^1\setminus[d_m,d_1] &\hbox{if $i=2m$}\\
               \{d_{{i+1\over 2}}\} & \hbox{if $i$ is odd}\\
               (d_{{i\over 2}},d_{{i\over 2}+1}) & \hbox{if $i$ is even.}
\endcases
$$

\begin{Cor}
$\P(W)^{ss}_{\la^k_d}=\P(W)^{ss}_{\la^{k^\p}_{d^\p}}$
if and only if there is an $i$ with $1\le i\le 2m$,
such that $I_i$ contains both $d/ k$ and $d^\p/ k^\p$.
\end{Cor}

We see that for the given action $\overline{\lambda}$
there are exactly $2m$ notions of stability. Denote
by $Q_i$, $i=1,...,2m$, the corresponding unpolarized GIT-quotients,
where $Q_{2m}=\emptyset$.
Then, for any $i=1,...,2m$, there is a $k$ with $Q_i=Q^k_2$.\par

\begin{Rem}
Bia\l ynicki-Birula and Sommese \cite{BS}
investigated $\C^*$-actions in a more
general context.
Specialized to our situation, their main result is the following:
 Let $\la$ be a $\C^*$-action on $W$ with
a decomposition of the dual space $W^\vee=\bigoplus_{i=1}^m W_i^\vee$
as above. The fixed point set of the induced $\C^*$-action on $\P(W)$ is
given by
$\bigcup_{i=1}^m \P(W_i)$. Set $F_i:=\P(W_i)$, and
define for each index $i$:
\begin{eqnarray*}
X_i^+&:=&\bigl\{\,x\in \P(W)\ \vert\ \mathop{\rm lim}_{z\lra 0} z\cdot x\in F_i
\,\bigr\}=\P(W_i\oplus\cdots\oplus W_m)\\
X_i^-&:=&\bigl\{\,x\in \P(W)\ \vert\ \lim_{z\lra \infty} z\cdot x\in
F_i\,\bigr\}=
\P(W_1\oplus\cdots \oplus W_i),
\end{eqnarray*}
and for $i\neq j$ set
$C_{ij}:=(X_i^+\backslash F_i)\cap (X_j^-\setminus F_j)$. This means
$C_{ij}$ is empty for $i\ge j$ and equal to
$\P(W_i\oplus\cdots\oplus W_j)\setminus (\P(W_i)\cup\P(W_j))$ for $i<j$.
We write $F_i<F_j$ when $C_{ij}\neq\emptyset$, i.e.
$$F_1<F_2<\cdots<F_m.$$
In the terminology of \cite{BS}, $F_1$ is the \sl source \rm and $F_m$ is
the \sl sink. \rm
For each $i$ with $1\le i\le m-1$,
one has a partition of $A:=\{1,...,m\}$:
$$A=A_i^-\cup A_i^+,\q \hbox{with } A_i^-:=\{1,...,i\} \hbox{ and } A_i^+=\{
i+1,...,m\},$$
and an associated
open set
$$U_i:=\bigcup_{\mu\in A_i^-, \nu\in A_i^+} C_{\mu\nu}.$$
The main theorem of \cite{BS} asserts that the $U_i$ are the only Zariski-open
$\C^*$-invariant subsets of $\P(W)$ not intersecting the fixed point set
whose quotients by the $\C^*$-action are compact.
One checks directly that $U_i$ is the set of
$\la^k_d$-semistable points for any pair $k,d$ with
$d/ k\in(d_i,d_{i+1}).$
\end{Rem}

\begin{Ex}
\label{C^*-ex}
Consider an action $\la$ of $\C^*$ on a finite dimensional vector space
$W$
such that the dual space decomposes as
 $W^\vee=W_1^\vee\oplus W_2^\vee$ with weights $d_1<d_2$.
If $d\in\Z$ and $k\in\Z_{>0}$ are such that $d_1<d/k <d_2$, then
the set of $\la^k_d$-semistable points is
$\P(W_1\oplus W_2)\backslash
\left(\P(W_1)\cup\P(W_2)\right)$ and the quotient $Q_{\la^k_d}$ is
naturally isomorphic to $\P(W_1)\times\P(W_2)$. The quotient map
$$\pi\colon\P(W_1\oplus W_2)\backslash
\left(\P(W_1)\cup\P(W_2)\right)
\subset \P(W_1\oplus W_2)\dasharrow \P(W_1)\times\P(W_2)$$
is the obvious one.
\end{Ex}

\begin{Claim}
The polarization induced by $\la^k_d$ on $\P(W_1)\times\P(W_2)$
is the equivalence class of the bundle
$\O_{\P(W_1)\times\P(W_2)}(kd_2-d,-kd_1+d)$.
In particular, for every $m,n\in \Z_{>0}$, the class
$[\O_{\P(W_1)\times\P(W_2)}(m,n)]$ occurs as an
induced polarization.
\end{Claim}

%\begin{pf}
{\it Proof}.
Let $L:=\O_{\P(W_1)\times\P(W_2)}(m,n)$ represent the induced polarization.
From the description of $\pi$ it follows that
$H^0(\pi^*L)^{{\la_d^k}}=\pi^*H^0(L)=S^mW_1\otimes S^nW_2$ is
the set of bihomogenous polynomials of bidegree $(m,n)$, for some  $m,n$.
If $S^mW_1\otimes S^nW_2$
occurs as an eigenspace of the induced $\C^*$-action on the space
$H^0(\O_{\P(W_1\oplus W_2)}(m\cdot n))$, then it must obviously  be an
eigenspace for
the character $\chi_{-(md_1+nd_2)+((m+n)/ k)  d}$.
Now invariance implies $md_1+nd_2-((m+n)/ k) d=0$,
which can be written as $m (kd_1-d) + n  (kd_2-d)=0$.
This yields the first assertion.\par
To prove the second part of the claim one has to find
positive integers $k$, $r$ and an integer $d$ such that the following
equations hold
\begin{eqnarray*}
kd_2-d &=& r  m\\
-kd_1+d &=& r  n;
\end{eqnarray*}
this results from a straightforward computation.

The other  quotients are $\P(W_1)$,  $\P(W_2)$ with
the obvious polarizations, and  $\emptyset$.
%\end{pf}

\subsection{Stability for actions of products of groups}

Consider now two reductive groups $G$, $H$ and a rational representation
$\rho\colon G\times H\lra \GL(W)$ in the finite dimensional space $W$.
We denote by $\gamma$ and $\la$ the induced representations of
$G$ and $H$, respectively.
Choose $n$ large enough in order to obtain an embedding
$j_n\colon Q_\gamma\hookrightarrow \P({S^nW}^G)$.
Since the actions of $G$ and $H$ commute, $\la$ induces actions of $H$
on $Q_\gamma$, on ${S^nW}^G$, and on
$\P({S^nW}^G)$; for these
actions $j_n$ is $H$-equivariant.
The action of $H$ on $Q_\gamma$ possesses a natural linearization
in
$j_n^*\O_{\P({S^nW}^G)}(1)$.
By \ref{BackGIT}(\ref{polind}), the corresponding concept of stability does
not depend
on the choice of $n$.
Let us denote the set of semistable points by $Q_\gamma^{ss}$ and the set
of polystable points by $Q_\gamma^{ps}$.

\begin{Prop}
\label{prodss}
The set of $\rho$-semistable points in the projective space $\P(W)$ is given by
$\P(W)^{ss}_\rho=\P(W)^{ss}_\gamma\cap \pi_\gamma^{-1}(Q^{ss}_\gamma)$,
and there exists a natural isomorphism
$Q_\gamma\catqot_\la H\cong Q_\rho$.
\end{Prop}

\begin{pf}
Suppose $x\in\P(W)$ is $\gamma$-semistable and its image
$\pi_\gamma(x)$ is $\la$-semistable in $Q_\gamma$.
If $n$ is large, $j_n(\pi_\gamma(x))$ is semistable in
$\P({S^nW}^G)$, so that there exists an integer $k\ge 1$
and a section $\overline{s}\in
H^0(\P({S^nW}^G),\O_{\P({S^nW}^G)}(k))^H$ not vanishing at \
$j_n(\pi_\gamma(x))$.
Identifying $\overline{s}\in {S^k({S^nW}^G)}^H$ with an element of
${S^{kn}W}^{G\times H}$, we obtain a $G\times H$-invariant section
in $\O_{\P(W)}(kn)$ not vanishing at $x$, hence $x$ is $\rho$-semistable.
\par
Conversely, suppose $x\in \P(W)_\rho^{ss}$. Then there exists, for some
$m\ge 1$, a section $s\in H^0(\P(W),\O_{\P(W)}(m))^{G\times H}$ with $s(x)\neq 0$.
Viewing $s\in {S^mW}^{G\times H}$ as an $H$-invariant element of ${S^mW}^G$,
we see that $x\in \P(W)^{ss}_\gamma\cap \pi_\gamma^{-1}(Q_\gamma^{ss}).$
This proves the first assertion.
\par
The second assertion follows immediately from the first one
and the universal property of the categorical quotient.
\end{pf}

The corresponding result for the polystable points is

\begin{Prop} The set of $\rho$-polystable points is
$\P(W)^{ps}_\rho=\P(W)^{ps}_\gamma\cap \pi_\gamma^{-1}(Q^{ps}_\gamma)$.
\end{Prop}

\begin{pf}
Let $x\in \P(W)$ be a $\gamma$-polystable point with $\pi_\gamma(x)\in
Q_\gamma^{ps}$.
By \ref{prodss}, $x$ is $\rho$-semistable. Choose
a $\rho$-polystable point
$y\in \overline{(G\times H)\cdot x}\cap \P(W)_\rho^{ss}$.
Projecting onto $Q_\gamma$, it follows that $\pi_\gamma(y)$ is contained
in $\overline{H\cdot \pi_\gamma(x)}$ and hence in $H\cdot \pi_\gamma(x)$,
because $\pi_\gamma(x)$ is polystable by assumption.
Therefore, there exists an $h\in H$ with
$\pi_\gamma(x)=h\cdot \pi_\gamma(y)=\pi_\gamma(h\cdot y)$.
But this means that the closures of the $G$-orbits of $x$ and $h\cdot y$
intersect, so that $G\cdot x\subset \overline{G\cdot (h\cdot y)}\cap
\P(W)^{ss}_\gamma$,
 since $x$ is $\gamma$-polystable.
In particular, $x\in \overline{(G\times H)\cdot y}\cap\P(W)_\rho^{ss}=
(G\times H)\cdot y.$ Hence $x$ is also $\rho$-polystable.
\par
To prove the converse, suppose $x$ is a $\rho$-polystable point.
We first show that $x$ is $\gamma$-polystable, too.
Let $y\in \overline{G\cdot x}\cap\P(W)_\gamma^{ss}$ be a $\gamma$-polystable
point. Since $\pi_\gamma(y)=\pi_\gamma(x)$, it follows from \ref{prodss} that
$\pi_\gamma(y)\in Q_\gamma^{ss}$.
Applying \ref{prodss} again, we see that $y\in\P(W)_\rho^{ss}$.
The orbit $(G\times H)\cdot x$ being closed in $\P(W)_\rho^{ss}$,
there exist $g\in G$  and $h\in H$ with $y=g\cdot h\cdot x$,
i.e.\ $x= h^{-1}\cdot g^{-1}\cdot y$.
Now $g^{-1}\cdot y$ is $\gamma$-polystable, hence $x$ is $\gamma$-polystable
too, because $\gamma$ and $\la$ commute.
Finally, we must show that $\pi_\gamma(x)\in Q_\gamma^{ps}$.
Choose $y$ such that $\pi_\gamma(y)\in\overline{H\cdot\pi_\gamma(x)}\cap
Q_\gamma^{ps}$.
 We may assume that $y$ is $\gamma$-polystable.
By what we have already proved, $y$ is $\rho$-polystable.
Now $\pi_\gamma(y)$ and $\pi_\gamma(x)$ are mapped to the
same point in $Q_\gamma\catqot_\la G=Q_\rho$.
But the projection $\pi_\rho\colon\P(W)_\rho^{ss}\lra Q_\rho$
separates closed $\rho$-orbits, thus $(G\times H)\cdot x=(G\times H)\cdot y$,
and therefore $H\cdot \pi_\gamma(x)=H\cdot \pi_\gamma(y)$ is closed
in $Q_\gamma^{ss}$.
\end{pf}

\subsection{Applications to $G\times\C^*$-actions}

Let $G$ be a reductive algebraic group possessing only the trivial
character, so that for any action of $G$ on a projective
variety $V$ and any
line bundle $L$ on $V$ there is at most one $L$-linearization
of the given action.
Consider a rational representation
$\rho$ of $G\times\C^*$ in the finite dimensional vector space $W$.
As above we denote by $\gamma$ and $\lambda$ the induced representations
of $G$ and $\C^*$, respectively, and  by $\overline{\rho}$,
$\overline{\gamma}$, and $\overline{\lambda}$ the induced action
of $G\times\C^*$, $G$, and $\C^*$ on $\P(W)$.
Let $\P(W)_i^s\subset\P(W)_i^{ps}\subset\P(W)_i^{ss}$
be the set stable, polystable, or semistable points w.r.t.\ the
$i$-th stability concept for the action $\overline{\lambda}$,
and let $I_i$, $i=1,...,2m$, be the associated intervals of rational
numbers.
The representation $\rho$ induces an action of $G$ on $Q^k_d$
which is equipped with a natural linearization in the ample line
bundle $L^k_d$, and there is no natural way to alter this linearization,
because $G$ does not possess a non-trivial character.
The corresponding concept of $G$-stability depends only on the
rational parameter $d/ k$.\par
Now fix a rational parameter $\eta:=d/ k\in I_i$ for some
index $i$.
A point $y\in Q_i$ is called \it
$\eta$-stable ($\eta$-polystable, $\eta$-semistable) \rm if it is
$G$-stable ($G$-polystable, $G$-semistable) w.r.t.\ the
$G$-linearized line bundle $L^k_d$ on $Q_i=Q_d^k$.
\par
Recall that every point $x\in\P(W)$ lies in $\P(W)_i^{ps}$
for a suitable index $i$; let $\pi_i(x)\in Q_i$ be its
image under $\pi_i\colon \P(W)_i^{ps}\lra Q_i$.

\begin{Thm}
\label{GITThm}
Fix a point $x\in\P(W)$. Then the following conditions are equivalent:
\par
{\rm i)} The point $x$ is $G$-semistable ($G$-polystable).\par
{\rm ii)} There exists an index $i$ and a parameter $\eta\in I_i$
such that $x\in \P(W)_i^{ss}$ ($x\in \P(W)_i^{ps}$) and
$\pi_i(x)$ is $\eta$-semistable ($\eta$-polystable).
\end{Thm}

\begin{pf} We explain the semistable case; the arguments in the
polystable case are similar.
Suppose first that $x\in \P(W)$ is $G$-semistable.
Choose $n$ large enough (cf.\ Section~\ref{BackGIT}) in order to obtain
a commutative diagram as in \ref{BackGIT}(\ref{eqemb}).
Since $\gamma$ and $\lambda$ commute, the representation
$\la^n\colon \C^*\lra \GL(S^nW)$
induces a representation
$\la^\p\colon \C^*\lra \GL({S^nW}^G).$
By \ref{C^*-ss2}, we find $k\in \Z_{>0}$ and $d\in\Z$
such that $\pi_\gamma(x)$ is semistable w.r.t.\ the stability
concept induced by $(\la^\p)^k_d$ on $Q_\gamma$.
Since $(\la^\p)^k_d$ is induced by the representation
$$\la^{nk}_d\colon \C^*\lra \GL\left(S^k(S^nW)\right),$$
we may replace $n$ by $kn$ and, therefore, assume that
$\pi_\gamma(x)$ is semistable w.r.t.\ the stability
concept induced by $(\la^\p)_d$ on $Q_\gamma$, for some integer $d$.
We now apply Proposition~\ref{prodss} to the representation
$$(\gamma^n\times \la^n_d)\colon G\times\C^*\lra \GL(S^nW).$$
(Note that this representation induces the action $\overline{\rho}$ on $\P(W)$.)
Since $x\in\P(W)$ is $\gamma$-semistable, it is also $\gamma^n$-semistable.
By construction, $\pi_\gamma(x)$ is semistable w.r.t.\ the induced
$\C^*$-action on $Q_\gamma$, and hence $x$ is $\gamma^n\times\la_d^n$-semistable
 by \ref{prodss}.
Applying \ref{prodss} the other way round, setting $\eta:=d/ k$
and choosing $i$ with $\eta\in I_i$,
it follows that $x\in\P(W)_i^{ss}$ and that $\pi_i(x)$ is $\eta$-semistable.
This settles the implication i)$\Rightarrow$ii).
\par
To prove the other implication suppose $x\in \P(W)$ fulfills
the assumptions of ii).
By definition and by Proposition~\ref{prodss}, this means that
there are $k\in\Z_{>0}$ and
$d\in\Z$ with $\eta=d/ k$
such that $x\in\P(W)$ is $\gamma^k\times\la^k_d$-semistable.
This implies that $x$ is $\gamma^k$- and hence $\gamma$-semistable.
This concludes the proof.
\end{pf}

\begin{Rem}
\label{ChainsofFlips}
At this point it becomes clear why chains of flips appear:
Let $G$, $\rho$, $\gamma$, and $\la$ be as above.
We have constructed a family
$(\gamma^k\times\la^k_d)$ of linearizations of the action $\overline{\rho}$
on $\P(W)$.
Each of these linearizations yields a GIT-quotient of $\P(W)$ by the action
$\overline{\rho}$.
This family of quotients can be constructed in another manner:
First take the $G$-quotient in order to obtain
a polarized variety $(\tilde{Q}:=\P(W)\catqot_{\gamma} G,[L])$.
The resulting $\C^*$-action on this variety yields a
family of quotients $Q_i$, $i=1,...,2n$, where $2n$ is usually (much) larger
than $2m$, the number of \sl unpolarized \rm $\C^*$-quotients of $\P(W)$
(see~\ref{FlipsEx}). But the family $Q_i$, $i=1,...,2n$,
coincides with the family $\P(W)\catqot_{\gamma^k\times \la_d^k}G\times\C^*$,
$k\in \Z_{>0}, d\in\Z$.
This phenomenon is responsible for the occurence of chains of flips
in these situations. It explains the question which was left open in
\cite{R}, 2.4 Remark (2), 2.5.
\end{Rem}

\begin{Ex}
\label{FlipsEx}
Let $W^\vee:=S^3{\C^2}^\vee\oplus {\C^2}^\vee$ and let $\SL_2(\C)$ act on
$W^\vee$ in the following way: Given $(f,p)\in W^\vee$ and $m\in \SL_2(\C)$,
we interpret $f$ and $p$ as functions on $\C^2$ and set
$(m\cdot f)(v):= f(m^t\cdot v)$ and $(m\cdot p)(v):=p(m^t\cdot v)$;
then we define $m\cdot (f,p) := (m\cdot f, m\cdot p).$
Let $\C^*$ act on $W^\vee$ by multiplication with $z^{d_1}$ on the first
factor and by multiplication with $z^{d_2}$ on the second one.
The quotient $V:=W^\vee\catqot\SL_2(\C)$ is of the form
$\Spec\C[I,J,D,R]$, where $I$, $J$, $D$, and $R$ are certain
bihomogenous polynomials of bidegrees $(2,2)$, $(3,3)$, $(4,0)$, and
$(1,3)$ in the coordinates of $S^3{\C^2}^\vee$
and ${\C^2}^\vee$.
Furthermore, $I$, $D$, and $R$ are algebraically independent,
and there is a relation
$$27 J^2={1\over 256}DR^2+I^3.$$
We examine the $\SL_2(\C)\times\C^*$-action on $\P(W)$.
The quotient $Q:=\P(W)\catqot\SL_2(\C)$ is given by $\Proj\C[I,J,D,R]$
where $I$, $J$, $D$, and $R$ have weights 4, 6, 4, and 4, respectively.
The ring $\C[I,J,D,R]_{(12)}$ is generated by its elements in degree 1, i.e.
by $I^3, I^2D, I^2R, ID^2, IR^2, IDR, J^2, D^3, D^2R, DR^2, R^3$;
hence there is an embedding $Q\hookrightarrow \P(S^{12}W^{\SL_2(\C)})$.
The $\C^*$-action on $Q$ can be extended to $\P(S^{12}W^{\SL_2(\C)})$ such
that the weights of the corresponding action on ${S^{12}W^{\SL_2(\C)}}^\vee$
are
$$6d_1+6d_2, 8d_1+4d_2, 5d_1+7d_2, 10 d_1+2d_2, 4d_1+8d_2,
7d_1+5d_2, 12 d_1, 9d_1+3d_2, 3d_1+9d_2.$$
For a point in $p\in Q$, $d_{\min}(p)$ and $d_{\max}(p)$ can take the
values $6d_1+6d_2$, $12d_1$, and $3d_1+9d_2$.
Hence, for $d_1\neq d_2$, there are
6 different notions of semistability on $Q$, hence 6 different
notions of
$\SL_2\times\C^*$-semistability on $\P(W)$, whereas there are only 4
different notions of
$\C^*$-semistability on $\P(W)$.

\end{Ex}

%%%%%%%%%%%%%%%%%%%%%%%%%%%%%%%%%%%%%%%%%%%%%%%%%%%%%%
\section{Oriented pairs and their moduli}

Let $X$ be a smooth projective variety over the field of complex numbers
and fix an ample divisor $H$ on $X$.
All degrees will be taken with respect to $H$ and
the corresponding line bundle will
be denoted by $\O_X(1)$.
Fix a torsion free
coherent sheaf $\E_0$ and a Hilbert polynomial $P$.
Finally, let $\Pic(X)$ be the Picard scheme of $X$
and choose a Poincar\' e line bundle $\L$ over $\Pic(X)\times X$.
If $S$ is a scheme and ${\frak E}_S$ a flat family of coherent sheaves
over $S\times X$, then there is a morphism $\det_S\colon S\lra \Pic(X)$
mapping a closed point $s$ to $[\det({\frak E}_{S\vert \{s\}\times X})]$.

We set $\L[{\frak E}_S]:=(\det_S\times\id)^*(\L)$; this line bundle depends only
 on the isomorphism class of the family ${\frak E}_S$.
The Hilbert polynomial of a sheaf $\F$ will be denoted by $P_\F$.
For any non-trivial
torsion free coherent sheaf $\F$ there is a unique subsheaf
$\F_{\max }$ for which $P_\F/\rk\F$ is maximal and whose rank is maximal among
the subsheaves $\F^\p$ with $P_{\F^\p}/\rk\F^\p$ maximal.
Set $\mu_{\max }(\F):=\mu(\F_{\max })$.

\subsection{Oriented pairs}

An \it oriented pair of type $(P,\L,\E_0)$ \rm is a triple $(\E,\eps,\phi)$
consisting of a torsion free coherent sheaf $\E$ with
Hilbert polynomial $P_\E=P$, a homomorphism $\eps\colon\det\E\lra\L[\E]$,
and a homomorphism $\phi\colon \E\lra \E_0$.
The homomorphisms $\eps$ and $\phi$ will be called
the \it orientation \rm and   the \it framing \rm of the pair $(\E,\eps,\phi)$.
Two oriented pairs $(\E_1,\eps_1,\phi_1)$ and $(\E_2,\eps_2,\phi_2)$ are
said to be \it equivalent\rm , if there is an isomorphism
$\Psi\colon\E_1\lra \E_2$ with $\eps_1
=\eps_2\circ\det\Psi$ and $\phi_1=\phi_2\circ\Psi$. When $\ker(\phi)\neq
0$, we set
$$\delta_{\E,\phi}:=P_\E-{\rk\E\over\rk\ker(\phi)_{\max}}P_{\ker(\phi)_{\max
}}.$$
An oriented pair $(\E,\eps,\phi)$ of type $(P,\L,\E_0)$ is \it semistable\rm ,
if either $\phi$ is injective, or $\eps$ is an isomorphism, $\ker(\phi)\neq 0$,
$\delta_{\E,\phi}\ge 0$, and for all non-trivial subsheaves $\F\subset\E$
$${P_\F\over \rk\F}-{\delta_{\E,\phi}\over\rk\F}\q\le\q {P_\E\over\rk\E}-
{\delta_{\E,\phi}\over\rk\E}.$$
\\
 The corresponding stability concept is slightly more
complicated: An oriented pair $(\E,\eps,\phi)$ of type $(P,\L,\E_0)$ is \it
stable\rm ,
if either $\phi$ is injective, or $\eps$ is an isomorphism, $\ker(\phi)\neq 0$,
$\delta_{\E,\phi}> 0$, and one of the following conditions holds:
\begin{enumerate}
\item For all non-trivial proper subsheaves $\F\subset \E$:
$${P_\F\over \rk\F}-{\delta_{\E,\phi}\over\rk\F}\q
<\q{P_\E\over\rk\E}-{\delta_{\E,\phi}\over\rk\E}.$$
\item $\phi\neq 0$, $\ker(\phi)_{\max}$ is stable, and
$\E\cong\ker(\phi)_{\max}\oplus\E^\p$, where the pair $(\E^\p,\phi)$ satisfies
\begin{eqnarray*}
{P_\F\over\rk\F}-{\delta_{\E,\phi}\over\rk\F}\q <\q{P_{\E^\p}\over\rk\E^\p}
-{\delta_{\E,\phi}\over\rk\E^\p} &\ & \hbox{$\forall$ proper subsheaves
$0\neq\F\subset\E^\p\ ,$}
\\
{P_\F\over\rk\F}\q <\q{P_{\E^\p}\over\rk\E^\p}-{\delta_{\E,\phi}
\over\rk\E^\p} &\ &
\hbox{$\forall$ proper subsheaves $0\neq \F\subset\E^\p\cap\ker(\phi)$.}
\end{eqnarray*}
\end{enumerate}
Our (semi)stability concept is related to the \sl parameter dependent \rm
(semi)stability concept of \cite{HL1} and \cite{HL2}
in the following way:
Let $\delta$ be a polynomial over the rationals with positive leading
coefficient.
 Recall that a pair $(\E,\phi)$ consisting of a torsion free coherent sheaf
$\E$ with $P_{\E}=P$ and a non-zero homomorphism $\phi\colon \E\lra\E_0$
is called \it (semi)stable w.r.t.\ $\delta$\rm , if for any non-trivial
proper subsheaf
$\F\subset\E$ the following conditions hold:
\begin{eqnarray*}
{P_\F\over\rk\F}-{\delta\over\rk\F}&(\le)&{P_\E\over\rk\E}-{\delta\over\rk\E
}\ ,\\
{P_\F\over\rk\F}&(\le)&{P_\E\over\rk\E}-{\delta\over\rk\E}
,\qquad\hbox{when $\F\subset\ker(\phi)$}.
\end{eqnarray*}
In this terminology, (semi)stable oriented pairs can be characterized as
follows:
\begin{Lem}
\label{HLcharac}
{\rm i)} An oriented pair $(\E,\eps,\phi)$ is semistable if and only if it
satisfies
 one of the following three conditions:
\begin{enumerate}
\item $\phi$ is injective.
\item $\E$ is semistable and $\eps$ is an isomorphism.
\item $\phi\neq 0$, $\eps$ is an isomorphism, and $(\E,\phi)$ is semistable
w.r.t.\ some $\delta>0$.
\end{enumerate}\par
{\rm ii)} An oriented pair $(\E,\eps,\phi)$ is stable if and only if it
satisfies
 one of the following four conditions:
\begin{enumerate}
\item $\phi$ is injective.
\item $\E$ is stable and $\eps$ is an isomorphism.
\item $\phi\neq 0$, $\eps$ is an isomorphism, and $(\E,\phi)$ is stable
w.r.t.\ some $\delta>0$.
\item $\phi\neq 0$, $\delta_{\E,\phi}>0$, $\eps$ is an isomorphism, and
$\E$ splits as $\ker(\phi)_{\max}\oplus\E^\p$, where $\ker(\phi)_{\max}$ is
stable and
$(\E^\p,\phi)$ is stable
 w.r.t.\ $\delta_{\E,\phi}$.
\end{enumerate}
\end{Lem}

We note that the stable oriented pairs appearing in
Lemma~\ref{HLcharac}.ii)4.
are precisely those pairs $(\E,\eps,\phi)$, for which
$\eps$ is isomorphic, $\phi\neq 0$, $\delta_{\E,\phi}>0$,
 the pair $(\E,\phi)$ is polystable w.r.t.\ $\delta_{\E,\phi}$,
and which have \sl only finitely many automorphisms\rm .
To see this, recall from
\cite{HL2} that for a given $\delta\in\Q[x]$, $\delta>0$, the polystable
 pairs
 $(\E,\phi)$
are those for which ${\cal E}$ splits in the form
$$\E\cong \E_1\oplus\cdots\oplus\E_{s-1}\oplus\E_s\ ,$$
where the sheaves $\E_1,...,\E_{s-1}$ are stable subsheaves of $\ker(\varphi)$,
$(\E_s,\phi)$ is a
stable pair w.r.t.\ $\delta$, and
$P_{\E_1}/\rk\E_1=\cdots=
P_{\E_{s-1}}/\rk\E_{s-1}=
P_{\E_s}/\rk\E_s-\delta/\rk\E_s$.

This makes our assertion obvious.

\begin{Rem}
\label{properlysemistablepairs}
Let $(\E,\eps,\phi)$ be a stable oriented pair of type 4.
(see~\ref{HLcharac}.ii)).
 Then $\delta_{\E,\phi}$
is the only rational polynomial with positive
leading coefficient w.r.t.\ which the pair $(\E,\phi)$ is semistable.
 This follows from
the equalities
\begin{eqnarray*}
{P_{\E^\p}\over\rk\E^\p}-{\delta_{\E,\phi}\over\rk\E^\p}&=
&{P_{\E}\over\rk\E}-{\delta_{\E,\phi}\over\rk\E}\ ,
\\
{P_{\ker(\phi)_{\max}}\over\rk\ker(\phi)_{\max}}&=
&{P_{\E}\over\rk\E}-{{\delta_{\E,\phi}}\over{\rk\E}}\stackrel{.}{}
\end{eqnarray*}

\end{Rem}
For all  stability concepts introduced so far, there are analogous notions of
\it slope-(semi)stability\rm . As usual, slope-stability implies stability and
semistability implies slope-semistability.
\par
Let $S$ be a noetherian scheme. \it A family of oriented pairs
parametrized by $S$ \rm is a quadruple $({\frak E}_S,\eps_S,
\widehat{\phi}_S,{\frak M}_S)$
consisting of a flat family ${\frak E}_S$ of torsion free
coherent sheaves
over the product $S\times X$, an invertible sheaf ${\frak M}_S$ on $S$,
a morphism
$\eps_S\colon \det{\frak E}_S\ra \L[{\frak E}_S]\otimes
\pi_S^*{\frak M}_S$,
and a morphism
$\widehat{\phi}_S\colon S^r{\frak E}_S\ra \pi_X^*S^r\E_0\otimes
\pi_S^*{\frak M}_S$
with $\widehat{\phi}_{S|\{s\}\times X}=S^r\phi_s$ for any
closed point $s\in S$ and a
suitable
$\phi_s\in\Hom({\frak E}_{S|\{s\}\times X},\E_0)$, so that
the pair
$(\eps_{S|\{s\}\times X},
\widehat{\phi}_{S|\{s\}\times X})$ is non-zero.

Two families $({\frak E}^i_S,\eps^i_S,\widehat{\phi}_S^i,{\frak M}_S^i)$,
$i=1,2$, are
called
\it equivalent\rm , if there exist an isomorphism
$\Psi_S\colon {\frak E}_S^1\lra {\frak E}_S^2$ and an isomorphism
${\frak m}\colon {\frak M}_S^1\lra {\frak M}_S^2$ such that
$(\id_{\L[{\frak E}_S^1]}\otimes \pi_S^*{\frak m})\circ
 \eps_S^1=\eps_S^2\circ \det\Psi$ and
$(\id_{\pi_X^*S^r\E_0}\otimes \pi_S^*{\frak m})\circ
\widehat{\phi}_S^1=\widehat{\phi}_S^2\circ S^r\Psi$.
\par
With these notions, we  define the functors
 $M^{ss}_{(P,\L,\E_0)}$ and $M^s_{(P,\L,\E_0)}$
of equivalence classes of families of semistable and stable oriented pairs
of type $(P,\L,\E_0)$.

\begin{Rem}
Though the definition of a family may appear a little odd at first sight,
it will become clear that families must be defined in this way
for technical reasons. Families of the above type are precisely those
 which are locally induced by the universal family on the
parameter space which we will construct in Section~\ref{ParSpace}.

The  functors defined above do depend on the
choice of the Poincar\' e bundle and  there is no natural way
to compare functors associated to different Poincar\' e bundles.
\end{Rem}
\subsection{A boundedness result}
\label{Bound}
Here we show that the family of isomorphism classes of torsion free
coherent sheaves occuring in oriented slope-semistable pairs of type
$(P,\L,\E_0)$  is
bounded. We use Maruyama's boundedness criterion:

\begin{Thm}\cite{Ma}
Let $C$ be some constant.
The set of isomorphism classes of torsion free coherent sheaves with
Hilbert polynomial $P$ and $\mu_{\max }\le C$ is bounded.
\end{Thm}

 \begin{Prop}
The set of isomorphism classes of torsion free  sheaves occuring in a
slope-semistable oriented pair of type $(P,\L,\E_0)$ is bounded.
\end{Prop}

\begin{pf}
Set $C:=\max\{\,\mu_{\max}(\E_0), \mu(\E)\,\}$.
Let $(\E,\eps,\phi)$ be a slope-semistable oriented pair of type
 $(P,\L,\E_0)$.
We claim that $\mu_{\max}(\E)\le C;$ in view of Maruyama's theorem, this
assertion proves
the proposition.
\par
Write a given non-trivial subsheaf $\F$ of $\E$ as an extension
$$0\lra \F\cap\ker(\phi)\lra \F\lra\phi(\F)\lra 0.$$
If $\F$ is entirely contained in the kernel of $\phi$, the definition of
slope-semistability
implies $\mu(\F)\le \mu(\E)\le C$.
If $\F$ is isomorphic to $\phi(\F)$, then obviously
$\mu(\F)\le\mu_{\max}(\E_0)\le C$.
In the remaining cases
\begin{eqnarray*}
\mu(\F) &=& {\mu(\F\cap\ker(\phi))\rk(\F\cap\ker(\phi))+
\mu(\phi(\F))\rk\phi(\F)\over\rk\F}\\
        &\le& {\rk(\F\cap \ker(\phi))\over\rk\F}\mu(\E)+
              {\rk\phi(\F)\over\rk\F}\mu_{\max}(\E_0) \le C.
\end{eqnarray*}
\end{pf}

\subsection{The parameter space for semistable oriented pairs}
\label{ParSpace}
By the boundedness result of the previous paragraph, there is a natural
number $m_0$ such that for all torsion free coherent sheaves
$\E$ occuring in a
semistable oriented pair, and for all $m\ge m_0$ the following
properties hold true: $\E(m)$ is globally generated and $H^i(X,\E(m))=0$
for $i>0$.
Let $V$ be a complex vector space of dimension $p:=P(m)$.
There exists a quasi-projective scheme ${\frak Q}$, the $\Quot$-scheme of
torsion free coherent quotient
sheaves of
$V\otimes\O_X(-m)$
with Hilbert polynomial $P$, and a universal
quotient on ${\frak Q}\times X$:
$$q_{\frak Q}\colon V\otimes\pi_X^*\O_X(-m) \lra {\frak E}_{\frak Q}.$$
Let $\n $ be the sheaf $\pi_{\frak Q *}(\det({\frak E}_{\frak Q})^\vee
\otimes \L[{\frak
E}_{\frak Q}])$. By the universal property of the Picard scheme, there is a line bundle
${\frak M}$ on ${\frak Q}$ such that
$$\det({\frak E}_{\frak Q})^\vee\otimes  \L[{\frak E}_{\frak Q}]\cong \pi_{\frak
Q}^*\frak M.$$
This implies that $\n$ is invertible and
$$\n\langle [q]\rangle\cong H^0(X,\det({\frak E}_{{\frak Q}\vert\{[q]\}\times
X}^\vee)\otimes\L[{\frak E}_{{\frak Q}\vert\{[q]\}\times X}])\ .$$
 Let ${\frak N}\buildrel
\over\lra {\frak Q}$ be the associated geometric line bundle.
The space ${\frak N}$ is a parameter space for equivalence classes
$[q\colon V\otimes\O_X(-m)\lra \E,\eps]$
consisting of a quotient \linebreak
$q\colon V\otimes\O_X(-m)\lra \E$
and an orientation
$\eps\colon \det(\E)\lra \L[\E]$.
Here two objects $(q_i\colon V\otimes \O_X(-m)\lra\E_i,\eps_i)$,
$i=1,2$, are
\it equivalent\rm , if there is an isomorphism $\Psi\colon \E_1\lra\E_2$
with $\Psi\circ q_1=q_2$ and $\eps_1=\eps_2\circ \det(\Psi)$.
\par
Next we have to construct a parameter space for all oriented
pairs.
We choose $m\ge m_0$ so large that $\E_0(m)$ is also globally generated.
Every oriented pair yields an element
in $K:=\Hom(V,H^0(\E_0(m)))$ and hence an element in $S^rK$.

On the projective bundle ${\frak P}:=\P(({\frak N}\times S^rK)^\vee)
\stackrel{\frak p}
{\lra}{\frak Q}$  there is a (nowhere vanishing) tautological section
$${\frak s}\colon
\O_{\frak P}\lra {\frak p}^*({\cal N}\oplus (S^rK\otimes{\cal O}_{\frak
Q}))\otimes\O_{\frak P}(1).$$
Let
$$q_{\frak P}\colon V\otimes \pi_X^*\O_X(-m)\lra {\frak E}_{\frak P}$$
be the pullback of the universal quotient on ${\frak Q}\times X$
to ${\frak P}\times X$.
We  view the pullback  $\pi_{\frak P}^*{\frak s}$ of ${\frak s}$ to
${\frak P}\times X$
as a pair consisting of a homomorphism
$$\eps_{\frak P}\colon \det({\frak E}_{\frak P})\lra \L[{\frak E}_{\frak
P}]\otimes
\pi_{\frak P}^*\O_{\frak P}(1)$$ and a homomorphism
$$\kappa_{\frak P}\colon S^rV\otimes\O_{{\frak P}\times X}\lra
 S^rH^0(\E_0(m))\otimes
\pi_{\frak P}^*\O_{\frak P}(1).$$

\begin{Rem}
\label{MorphtoN}
For a scheme $S$, giving a morphism
$f\colon S\lra {\frak P}$ is equivalent to giving a map
 $\overline{f}\colon S\lra
{\frak Q}$ -  which yields the family ${\frak
E}_S:=(\overline{f}\times\id_X)^*{\frak
E}_{\frak Q}$ - , a line bundle ${\frak M}_S$ on $S$, and
homomorphisms
$$\eps_S\colon \det({\frak E}_S)\lra \L[{\frak E}_S]\otimes
\pi_S^*{\frak
M}_S\ ,$$
$$\kappa_S\colon S^rV\otimes\O_{S\times X}\lra
S^rH^0(\E_0(m))\otimes\pi_S^*{\frak
M}_S\ $$
on $S\times X$ such that the pair $(\eps_{S|\{s\}\times
X},\kappa_{S|\{s\}\times X})$ is
non-zero for every closed point $s\in S$. Of course, for the morphism $f$
determined by
$\overline{f}$ and $(\eps_S,\kappa_S,{\frak M}_S)$,  we have
$\overline{f}={\frak p}\circ
f$, and there is an isomorphism
${\frak m}\colon {\frak M}_S\lra \overline{f}^*\O_{\frak P}(1)$ such that
$$({\rm id}_{\L[{\frak E}_S]}\otimes \pi_S^*{\frak
m})\circ\eps_S=(f\times{\rm id}_X)^*(\eps_{{\frak P}})\ ,$$
$$ ({\rm id}_{\pi_X^*S^rH^0(\E_0(m))}\otimes\pi_S^*{\frak m})\circ
\kappa_S=(f\times{\rm id}_X)^*(\kappa_{\frak P})\ .$$
\end{Rem}
Our parameter space ${\frak B}$ will be a closed subscheme of
${\frak P}$ whose closed points are of the form
$[[q\colon V\otimes\O_X(-m)\lra \E,\eps], S^rk]$, with
$[q,\eps]\in {\frak N}$ and $k\in K$, such that there is a map
$\phi\colon \E\lra \E_0$ making the following diagramm
commutative:

%%%%%%%%%%%%%%%%%%%%%%%%%%%%%%%%%%%%%
%\begin{equation}
%\label{eqemb}
%\begin{array}{c}
\begin{center}
\unitlength=1mm
\begin{picture}(70,24)(0,8)
\put(0,29){$V\otimes{\cal O}_X(-m)$}
\put(28,30){\vector(1,0){30}}
\put(42,32){${\scriptstyle q}$}
\put(62,29){${\cal E}$}
\put(7,25){\vector(0,-1){10}}
\put(63,25){\vector(0,-1){10}}
\put(-15,9){$H^0({\cal E}_0(m))\otimes{\cal O}_X(-m)$}
\put(28,10){\vector(1,0){30}}
\put(62,9){${\cal E}_0$}
\put(2,20){$\scriptstyle k$}
\put(42,12){$\scriptstyle ev$}
\put(66,20){${\scriptstyle \varphi}$}
\end{picture}
%
%\end{array}$
%\end{equation}
\end{center}
%%%%%%%%%%%%%%%%%%%%%%%%%%%%%%%%%%%%%%%%

Scheme-theoretically, ${\frak B}$ is constructed as follows:
On ${\frak P}\times X$, there is a homomorphism
$$\overline{\phi}_{\frak P}\colon S^rV\otimes\pi_X^*\O_X(-rm)\lra
\pi_X^*S^r\E_0\otimes\pi_{\frak P}^*\O_{\frak P}(1).$$
Set $\widehat{\cal G}:=\ker(S^rq_{\frak P})$, choose $n\ge m$
large enough so that
$\widehat{\cal G}_{\vert \{b\}\times X}(n)$ is globally generated
and without higher cohomology for any closed point $b\in {\frak P}$,
and let
$$\widehat{\gamma}\colon {\cal G}:=\widehat{\cal G}\otimes\pi_X^*
\O_X(n)\lra \pi_X^*
S^r\E_0(n)\otimes\pi_{\frak P}^*\O_{\frak P}(1)$$ be the induced
homomorphism.
We first define a scheme $\widehat{\frak B}$ whose closed points
are those elements $b\in {\frak P}$ for which
$\widehat{\gamma}_{|\{b\}\times X}$
is the zero map.
Since ${\cal G}_{|\{b\}\times X}$ and $S^r\E_0(n)$ are globally generated
for any closed point $b\in {\frak P}$, the scheme $\widehat{\frak B}$
is the zero locus of the following homomorphism between locally free
sheaves:
$$\gamma:=\pi_{{\frak P}*}(\widehat{\gamma})\colon
\pi_{{\frak P}*}{\cal G}\lra \pi_{{\frak P}*}(\pi_X^*S^r\E_0(n)\otimes
 \pi_{\frak
P}^*\O_{\frak P}(1))=H^0(S^r\E_0(n))\otimes \O_{\frak P}(1).$$ The
scheme ${\frak B}$ we
are looking for is the scheme-theoretic intersection of
$\widehat{\frak B}$ with the
image in ${\frak P}$ of the weighted projective  bundle  associated with
the vector
bundle ${\frak N}\times K$ over ${\frak Q}$.  There exists a universal family
$({\frak
E}_{\frak B},\eps_{\frak B},\widehat{\phi}_{\frak B},{\frak M}_{\frak B})$:
${\frak
M}_{\frak B}$ is the restriction of
$\O_{\frak P}(1)$ to ${\frak B}$,
$q_{\frak B}$ and $\eps_{\frak B}$ are the restrictions of $q_{\frak P}$ and
$\eps_{\frak P}$, and $\widehat{\phi}_{\frak B}$ is induced by the restriction
of $\widehat{\phi}_{\frak P}$ which factorizes through
$S^r{\frak E}_{\frak B}$  by definition.
In the following, a closed point $b=[[q\colon V\otimes\lra \O_X(-m),\eps],
S^rk]\in
{\frak B}$ will be denoted by
$[q,\eps,\phi]$; here $\phi$ is the unique framing on $\E$ induced by $k$.

\begin{Rem}
By construction, a morphism $\widehat{f}\colon S\lra {\frak P}$ factorizes
through ${\frak B}$ if and only if it factorizes through the image
of the associated weighted projective bundle of ${\frak N}\times K$, and
$(\widehat{f}\times\id_X)^*(\widehat{\phi}_{\frak P})$ is identically
 zero on the kernel
of the map $(\widehat{f}\times \id_X)^*(S^rq_{\frak P})$.
\end{Rem}

On the parameter space ${\frak B}$, there is a natural action (from the
right) of
the group $\SL(V)$.
To define this action, it suffices to construct a $\SL(V)$-action on
${\frak P}$ which leaves the scheme ${\frak B}$ invariant.
The standard representation of $\SL(V)$ on $V$ gives us
 the homomorphism
$$\Gamma\colon V\otimes \O_{{\frak Q}\times\SL(V)\times X}\lra V\otimes
\O_{{\frak
Q}\times\SL(V)\times X}.$$ Moreover, on ${\frak Q}\times\SL(V)\times X$
there is the
pullback of the universal quotient
$$\pi_{{\frak Q}\times X}^*(q_{\frak Q})\colon V\otimes
 \pi_X^*\O_X(-m)\lra
\pi_{{\frak Q}\times X}^*{\frak E}_{\frak Q}.$$
By the universal property of the $\Quot$-scheme,
$\pi_{{\frak Q}\times X}^* (q_{\frak
Q})\circ \bigl(\Gamma\otimes\id_{\pi_X^*\O_X(-m)}\bigr)$ yields a morphism
$\overline{f}\colon {\frak Q}\times\SL(V)\lra {\frak Q}$ such that there is
a well-defined isomorphism
$$\Psi_{{\frak Q}\times \SL(V)}\colon (\overline{f}
\times\id{}_X)^*  {\frak E}_{\frak
Q}\lra \pi^*_{{\frak Q}\times X}{\frak E}_{\frak Q}$$ with $\Psi_{{\frak
Q}\times\SL(V)}\circ (\overline{f}\times\id_X)^*(q_{\frak Q})=
\pi_{{\frak Q}\times
X}^*(q_{\frak Q})\circ
\bigl(\Gamma\otimes\id_{\pi_X^*\O_X(-m)}\bigr)$. Let
$\Psi_{{\frak P}\times \SL(V)}$  be
the pullback of
$\Psi_{{\frak Q}\times\SL(V)}$ to
${\frak P}\times \SL(V)\times X$, and  set ${\frak M}_{{\frak P}\times\SL(V)}:=
\pi_{\frak P}^*\O_{\frak P}(1)$,
\begin{eqnarray*}\eps_{{\frak P}\times\SL(V)}&:= &\pi_{{\frak P}\times X}^*
(\eps_{{\frak
P}})\circ \det\Psi_{{\frak P}\times\SL(V)}\ ,\\
\kappa_{{\frak P}\times\SL(V)} &:=& \pi^*_{{\frak P}\times X}
(\kappa_{\frak P})\circ
S^r\left(({\frak p}\times\id{}_{\SL(V)\times X})^*\Gamma\right)\ .
\end{eqnarray*}
By Remark~\ref{MorphtoN}, the data $\overline{f}$ and ($\eps_{{\frak
P}\times\SL(V)},
\kappa_{{\frak P}\times \SL(V)}, {\frak M}_{{\frak P}\times\SL(V)})$ define
an action
$$f\colon {\frak P}\times \SL(V)\lra {\frak P}.$$

\begin{Prop}
\label{LocUnivProp}
Let $S$ be a noetherian scheme and let $({\frak E}_S,\eps_S,\widehat{\phi}_S,
 {\frak
M}_S)$ be a family of semistable oriented pairs parametrized by $S$.
Then $S$ can be covered by open subschemes $S_i$ for which there exist
morphisms $\beta_i\colon S_i\lra {\frak B}$
such that the restricted families
$({\frak E}_{S|S_i},\eps_{S|S_i},\widehat{\phi}_{S|S_i},{\frak M}_{S|S_i})$
are equivalent
to the
pullbacks of $({\frak E}_{\frak B}, \eps_{\frak B}, \widehat{\phi}_{\frak B},
 {\frak
M}_{\frak B})$ via the maps $\beta_i\times\id_X$.
\end{Prop}
\begin{pf}
The scheme $S$ can be covered by open subschemes $S_i$ such that
the family ${\frak E}_{S\vert S_i}$ over $S_i\times X$ can be written as a
family
of quotients:
$$q_{S_i}\colon V\otimes \pi_X^*\O_X(-m)\lra {\frak E}_{S|S_i}.$$
Each $q_{S_i}$ defines a morphism $\overline{f}_i\colon S_i\lra {\frak Q}$
such that there is a well defined isomorphism
$\Psi_{S_i}\colon {\frak E}_{S_i}:=(\overline{f}_i\times\id_X)^*
{\frak E}_{\frak Q}\lra {\frak E}_{S|S_i}$.
Define ${\frak M}_{S_i}:={\frak M}_{S|S_i}$,
$$\begin{array}{cl}
\eps_{S_i}\colon &  \det({\frak E}_{S_i})\textmap{\det\Psi_{S_i}}
\det{\frak E}_{S|S_i}\textmap{\eps_{S|S_i}} \L[{\frak E}_{S\vert S_i}]
 \otimes
\pi_{S_i}^*{\frak M}_{S_i}\ \hbox{,}\\
\widehat{\phi}_{S_i}\colon &   S^r{\frak E}_{S_i}\textmap{S^r\Psi_{S_i}}
S^r{\frak
E}_{S|S_i}\textmap{\widehat{\phi}_{S|S_i}}\pi_X^*S^r\E_0\otimes\pi_{S_i}^*
 {\frak
M}_{S_i}.
\end{array}$$
The homomorphism $\widehat{\phi}_{S_i}$
yields a homomorphism
$$\overline{\kappa}_{S_i}\colon S^rV\otimes\O_{S_i\times X}\lra
\pi_X^*S^r\E_0(m)\otimes\pi_{S_i}^*{\frak M}_{S_i}$$
and hence a homomorphism
$$\kappa_{S_i}:=\pi_{S_i}^*\pi_{S_i*}(\overline{\kappa}_{S_i})\colon S^rV
\otimes\O_{S_i\times X}
 \lra S^rH^0(\E_0(m))\otimes \pi^*_{S_i}{\frak M}_{S_i};$$
here we have used the fact that our definition of a family implies that the map
$$\pi_{S_i*}(\overline{\kappa}_{S_i})\colon S^rV\otimes\O_{S_i}\lra
H^0(S^r\E_0(m))\otimes{\frak M}_{S_i}$$
  factorizes through $S^rH^0(\E_0(m))\otimes {\frak
M}_{S_i}$.

 By Remark~\ref{MorphtoN}, the quadruple
$(\overline{f}_i,\eps_{S_i},\kappa_{S_i},{\frak M}_{S_i})$ determines a morphism
$\beta_i\colon S_i\lra {\frak P}$.
It is clear that the morphism
$\beta_i$ factorizes through ${\frak B}$ and that the family
$({\frak E}_{S_i},\eps_{S_i}, \widehat{\phi}_{S_i}, {\frak M}_{S_i})$ is the
 pullback of
the universal family by $\beta_i\times \id_X$.
The family $({\frak E}_{S_i},\eps_{S_i}, \widehat{\phi}_{S_i},{\frak M}_{S_i})$
is  equivalent to $({\frak E}_{S|S_i},\eps_{S|S_i},\widehat{\phi}_{S|S_i},
{\frak M}_{S|S_i})$ by construction.
\end{pf}
Let ${\frak B}^{\mathop{\rm iso}}$ be the open subscheme of oriented pairs
 $[q,\eps,\phi]$ for which
$H^0(q(m))$ is an isomorphism. The maps constructed in
the above proof factorize through ${\frak B}^{\mathop{\rm iso}}$.

\begin{Prop}
\label{GlueTog}
Let $S$ be a noetherian scheme and let $\beta_i\colon S\lra {\frak
B}^{\mathop{\rm iso}}$,
 $i=1,2$,
be two morphisms such that the pullbacks of $({\frak E}_{\frak B},
\eps_{\frak B},
\widehat{\phi}_{\frak B},{\frak M}_{\frak B})$ via the maps
$(\beta_i\times\id_X)$
are equivalent families. Then there exists an \'etale cover $\eta\colon
T\lra S$ and a
morphism $g\colon T\lra \SL(V)$ such that
$\beta_1\circ\eta=(\beta_2\circ\eta)\cdot g.$
\end{Prop}

\begin{pf}
Denote the two families by $({\frak E}_S^i,\eps_S^i,\widehat{\phi}_S^i,
{\frak M}_S^i)$,
and let $\Psi_S\colon {\frak E}_S^1\lra {\frak E}_S^2$ be the corresponding
 isomorphism.
The bundles ${\frak E}_S^i$  can be written as quotients
$q_S^i\colon V\otimes \pi_X^*\O_X(-m)\lra {\frak E}_S^i$, and there is a
morphism $g_S\colon S\lra \GL(V)$ making the following diagramm commutative:
%%%%%%%%%%%%%%%%%%%%%%%%%%%%%%%%%%%%%%%%%%
\begin{center}
\unitlength=1mm

\begin{picture}(70,24)(6,9)
\put(0,29){$V\otimes\pi_X^*{\cal O}_X(-m)$}
\put(32,30){\vector(1,0){20}}
\put(37,32){${\scriptstyle g_S\otimes{\rm id}}$}
\put(56,29){$V\otimes\pi_X^*{\cal O}_X(-m)$}
\put(7,25){\vector(0,-1){10}}
\put(64,25){\vector(0,-1){10}}
\put(5,9){${\frak E}^1_S$}
\put(16,10){\vector(1,0){42}}
\put(62,9){${\frak E}^2_S$}
\put(2,20){$\scriptstyle q^1_S$}
\put(35,12){$\scriptstyle \Psi_S$}
\put(66,20){${\scriptstyle q^2_S}$}
\end{picture}
\end{center}
%%%%%%%%%%%%%%%%%%%%%%%%%%%%%%%%%%%%%%%%%
As in the proof of \cite{HL1}, Lemma 1.15, one constructs an
\'etale cover $\eta\colon
T\lra S$ such that there is a morphism ${\frak d}\colon T\lra \C^*$ with
$({\frak
d}(t))^p=\det(g_S(\eta(t)))$ for any closed point $t\in T$. Now
define $g:={\frak d}\cdot
(g_S\circ \eta)$.  In view of the description of the $\SL(V)$-action at the
 beginning of
this section, the assertion is obvious.
\end{pf}

\subsection{The GIT-construction}
\label{GITconst}
Let ${\frak A}$ be the union of the finitely many components of $\Pic(X)$
meeting the image of $\det_{\frak B}\colon {\frak B}\lra \Pic(X)$.
We may choose $m$ so large that the restriction of the line bundle
$\L_{\vert{\frak
A}\times X}\otimes\pi_X^*\O_X(rm)$ to $\{a\}\times X$ is globally generated
and without
higher cohomology for any closed point $a\in {\frak A}$. The direct image sheaf
$\pi_{{\frak A}_*}(\L_{\vert{\frak A}\times X}\otimes\pi_X^*\O_X(rm))$ is then
locally
free and commutes with base change. The same holds for ${\cal H}om(\bigwedge^r
V\otimes\O_{\frak A},\pi_{{\frak A}*}(\L_{\vert{\frak A}\times X}
\otimes\pi_X^*\O_X(rm)))$; let ${\frak H}$ be the geometric
vector bundle associated to this locally
free sheaf.
Consider the homomorphism
$$\sigma_{\frak N}\colon \bigwedge^rV\otimes
 \O_{{\frak N}\times X}\lra
 \det{\frak E}_{\frak N}\otimes
\pi_X^*\O_X(rm)\stackrel{\eps_{\frak N}}{\lra}  \L[{\frak E}_{\frak N}]
\otimes\pi_X^*\O_X(rm).$$
By the
universal property of the scheme ${\frak H}$, the pushforward
 $\pi_{{\frak N}*}
(\sigma_{\frak N})$ determines a morphism of schemes
${\frak N}\lra {\frak H}$ and hence a
morphism
${\frak N}\times S^rK\lra {\frak H}\times S^rK$. Let ${\frak Z}$ be the
vector bundle
$({\frak H}\times S^rK)^\vee$ over
${\frak A}$, and denote by $\P({\frak Z})$ the associated projective bundle.
 $\P({\frak
Z})$ can be polarized by tensorizing ${\cal O}_{\P({\frak Z})}(1)$ with the
pull back
of a very ample line bundle from  ${\frak A}$.

On $\P({\frak Z})$ there is a natural action  of the group $\SL(V)$ from
the right,
which   is trivial on the base ${\frak A}$
and  admits a canonical linearization in the polarizing line bundle.
We have a natural morphism
$$\iota\colon {\frak B}\hookrightarrow {\frak P}\lra
\P({\frak Z}) $$
which is equivariant w.r.t.\ the given actions.

Let us describe the effect of  $\iota$
on closed points:
Given $b\in {\frak B}$, let $(\E_b,\eps_b,\phi_b)$ be the oriented pair
induced by the restriction of $({\frak E}_{\frak B},\eps_{\frak B},
\widehat{\phi}_{\frak B})$ to $\{b\}\times X$, i.e., $\E_b$ and $\eps_b$ are the
restrictions of ${\frak E}_{\frak B}$ and $\eps_{\frak B}$ and
$\phi_b$ is
 a framing with
$S^r\phi_b=\widehat{\phi}_{{\frak B}\vert \{b\}\times X}$ ($\phi_b$ is
unique up to an
$r$-th root of unity). The point $b$ is
 mapped to $[\L[\E_b], h,S^rk]$ with
$$h\colon \bigwedge^rV\lra H^0(\det(\E_b)(rm))\textmap{H^0(\eps_b(rm))}
H^0(\L[\E_b](rm))$$
 and $k=H^0((\phi_b\circ q)(m))$.
A point in $\P({\frak Z})$  is $\SL(V)$-\it (semi)stable \rm if
it is semistable in the projective space
$\P((\Hom(\bigwedge^rV,H^0(\L_{\vert\{a\}\times X}(rm)))\oplus
S^rK)^\vee),$ where
 $a$ is its image in ${\frak A}$.
\par
Let ${\frak B}^{ss}$ (${\frak B}^s$)
be the open subscheme of points
$[q,\eps,\phi]$ such that the triple
$(\E,\eps,\phi)$ is a semistable (stable)
oriented pair and such that  the homomorphism $H^0(q(m))\colon V\lra
H^0(\E(m))$ is
an isomorphism.
\par
\begin{Thm}
\label{MyStabCrit}
For $m$ large enough,
${\frak B}^{ss}=\iota^{-1}(\P({\frak Z})^{ss})$,
 and ${\frak B}^s=\iota^{-1}(\P({\frak Z})^{s})$.
\end{Thm}

Before we can start with the proof, we have to recall some definitions and
results from \cite{HL1} and \cite{HL2}.
Let $(\E,\phi)$ be a pair consisting of a torsion free coherent sheaf $\E$
with $P_\E=P$
and
a non-trivial framing $\phi$.
\par
Let $\overline{\delta}$ be any positive rational number.
The pair $(\E,\phi)$ is called \it sectional (semi)stable w.r.t.\ \rm
$\overline{\delta}$, if there is a
subspace $V\subset H^0(\E)$ of dimension $\chi(\E)=P(0)$ such that the
following conditions are satisfied:
\begin{enumerate}
\item For all non-trivial submodules $\F$ of $\ker(\phi)$:
$$(\rk\E)\dim \left(H^0(\F)\cap V\right)(\le)
\rk\F(\chi(\E)-\overline{\delta}).$$
\item For all non-trivial submodules $\F\neq \E$:
$$(\rk\E)\dim \left(H^0(\F)\cap V\right)(\le) \rk\F(\chi(\E)-\overline{\delta})+
(\rk\E)\overline{\delta}.$$
\end{enumerate}
\par
Then we have the following result \cite{HL2}{, Th.\ 2.1}:

\begin{Thm}
\label{SecStab}
For any polynomial $\delta$, there exists a natural number $m_1$
 such that for all $m\ge m_1$
the following conditions are equivalent for a pair $(\E,\phi)$:
\par
{\rm i)} $(\E,\phi)$ is (semi)stable w.r.t.\ the polynomial $\delta$.\par
{\rm ii)} $(\E,\phi)(m)$ is sectional (semi)stable w.r.t.\ $\delta(m)$.
\end{Thm}

Let $(q\colon V\otimes\O_X(-m)\lra\E,\phi)$ be a pair consisting of
a \sl generically \rm surjective map $q$ of $V\otimes \O_X(-m)$ to a
torsion free sheaf $\E$ with $P_\E=P$ and a non-zero homomorphism
$\phi\colon\E\lra\E_0$.
We can associate to this pair an element
$([h], [k])\in\P(H^\vee)\times\P(K^\vee)$, where
$H:=\Hom(\bigwedge^rV,H^0(\L[\E](rm)))$.
There is a natural $\SL(V)$-action on
$\P(H^\vee)\times\P(K^\vee)$ which
can be linearized in every sheaf $\O(a_1,a_2)$,
where $a_1$ and $a_2$ are positive integers.
Define $\nu:=a_2/a_1$ and $\overline{\delta}:= p\nu/(\rk\E +\nu)$.
The proof of \cite{HL1}, Proposition 1.18 is valid in any dimension
 and yields
the following

\begin{Thm}
\label{StabCrit}
Let $(q\colon V\otimes \O_X(-m)\lra\E,\phi)$ be as above.
The associated element $([h], [k])$ is (semi)stable
w.r.t.\ the linearization in $\O(a_1,a_2)$ if and only if the
following two conditions are satisfied:
\par
{\rm i)} The homomorphism $H^0(q(m))$ is injective.\par
{\rm ii)} The pair $(\E,\phi)(m)$ is sectional (semi)stable w.r.t.\
$\overline{\delta}$.
\end{Thm}

We also need the following obvious observation:
\begin{Lem}
\label{EasyStabCrit}
Let $(q\colon V\otimes\O_X(-m)\lra \E,\phi)$ be as above.
The following conditions are
equivalent:
\par
{\rm i)} The homomorphism $k=H^0((\phi\circ q)(m))$ is injective.\par
{\rm ii)} The associated element $[k]\in \P(K^\vee)$ is stable.
\end{Lem}

After these preparations, we return to our situation.
Let ${\frak B}_0\subset {\frak B}$ be
the open set of all oriented pairs $[q,\eps,\phi]$
for which $\E$ is semistable, and define for each polynomial $\delta$
the set ${\frak B}_\delta$ as the open set of oriented pairs $(\E,\eps,\phi)$
with $\phi\neq 0$ such that $(\E,\phi)$ is semistable w.r.t.\ $\delta$.
The union ${\frak B}^\p:={\frak B}_0\cup\bigcup {\frak B}_\delta$
is quasi-projective,
hence quasi-compact, so that there exist finitely many polynomials, say,
$\delta_1$,...,$\delta_s$ with
${\frak B}^\p={\frak B}_0\cup {\frak B}_{\delta_1}\cup\cdots\cup
{\frak B}_{\delta_s}$.
Let $M$ be some constant.
By \cite{Ma}{, Theorem 1.7}, the set of points $b \in {\frak B}$
such that $\mu_{\max}({\frak E}_{{\frak B}\vert\{b\}\times X})\le M$
is open. Since ${\frak B}$ is quasi-compact,
there is a constant $\mu_0$ such that
$\mu_{\max}({\frak E}_{{\frak B}\vert\{b\}\times X})\le \mu_0$ for all
$b\in {\frak B}$.
We also know that the family ${\frak Ker}$ of kernels of framings of semistable
 oriented pairs
is bounded. It follows that $\mu_{\max}(\ker(\phi))$, for $\ker(\phi)\in
{\frak Ker}$,
can only take finitely many values.
As in \cite{HL2}, Lemma 2.7, this implies that there are only finitely many
polynomials of the form $P_{\ker(\phi)_{\max}}$.
In particular, there are only finitely many polynomials of the form
$$P_\E-(\rk\E/\rk\ker(\phi)_{\max})P_{\ker(\phi)_{\max}}.$$
We assume in the following that these polynomials are among
$\delta_1,...,\delta_s$, and
 that the chosen $m$ is large enough, so that  Theorem~\ref{SecStab} holds
for all
 $\delta_i$ and
set $\overline{\delta}_i:=\delta_i(m)$.

\begin{Thm}
\label{TheProp}
Suppose $m$ is sufficiently large. Let $[q,\eps,\phi]\in {\frak B}$
be a pair with $\phi\neq 0$
which is not (semi)stable.
Then there is no positive rational number
$\overline{\delta}$ such that $(\E,\phi)(m)$ is sectional (semi)stable
w.r.t.\ $\overline{\delta}$.
 \end{Thm}

\begin{pf}
Denote by ${\frak S}$ the bounded set of equivalence classes of pairs
$(\E,\phi)$
for which there is an element $[q,\eps,\phi]\in {\frak B}$.
\par
By the above, any pair $(\E,\phi)\in {\frak S}$ satisfies
$\mu_{\max}(\E)\le\mu_0$.
Let $\tilde\delta$ be a rational polynomial of degree $\dim X-1$
whose leading coefficient $\tilde{\delta}_0$ satisfies
$\mu(\E)+\tilde{\delta}_0\ge \max\{\,0,\mu_0\,\}$.
One can now copy the proof of \cite{HL2}{, page 305}, to show that
there is a constant $C$ such that for any submodule $(\tilde{\E},\tilde{\phi})$
of a pair $(\E, \phi)\in {\frak S}$ either
$\vert \deg(\tilde{\E})-\rk\tilde{\E}\mu(\E)\vert <C$,
or for all $m$ large enough
\begin{eqnarray*}
{h^0(\tilde{\E}(m))\over \rk\tilde{\E}}-
{\tilde{\delta}(m)\over \rk\tilde{\E}}&<&
 {P_\E(m)\over \rk\E}-
{\tilde{\delta}(m)\over \rk\E}\qquad
\hbox{if $\tilde{\E}\not\subset \ker(\phi)\ ,$}\\
\\
{h^0(\tilde{\E}(m))\over \rk\tilde{\E}} &<&{P_\E(m)\over \rk\E}-
{\tilde{\delta}(m)\over
 \rk\E}\qquad\hbox{otherwise.}\\
\end{eqnarray*}
Recall that a submodule $\tilde{\E}\subset \E$ is called \it saturated\rm ,
if the quotient $\E/\tilde{\E}$ is torsion free.
The family of saturated submodules $\tilde{\E}$ of modules $\E$
with $(\E,\phi)\in {\frak S}$
satisfying $\vert \deg(\tilde{\E})-\rk\tilde{\E}\mu(\E)\vert <C$ is bounded
(\cite{HL2}{, Lemma 2.7}).
Denote this family by $\tilde{\frak S}$.
There are only finitely many
possibilities for the Hilbert polynomials of those submodules.
Let
$\delta_j^\p$ be the finite family of polynomials of the form
$P_\E-(\rk\E/ \rk\E^\p)P_{\E^\p}$ where $\E^\p$ is a saturated submodule
of $\ker(\phi)$ for some $(\E,\phi)\in \tilde{\frak S}$, and
$\delta_k^{\p\p}$ be the finite family of polynomials of the form
$(\rk\E^{\p\p} P_\E- \rk\E P_{\E^{\p\p}})/(\rk\E-\rk\E^{\p\p})$
where $\E^{\p\p}$ is a saturated submodule of a pair $(\E,\phi)\in
\tilde{\frak S}$
not contained in the kernel of $\phi$.
We may assume that $\tilde{\delta}$, the $\delta_j^\p$'s and the
$\delta_k^{\p\p}$'s
with positive leading coefficients are among $\delta_1,...,\delta_s$.
Next, we choose $m$ large enough, so that
$\tilde{\E}(m)$ is globally generated and has no higher cohomology for all
 $\tilde{\E}\in \tilde{\frak S}$.
Let $(\E,\phi)$ be a pair which is not semistable
w.r.t.\ any of the polynomials $\delta_1,...,\delta_s$.
This is equivalent to $(\E,\phi)(m)$ not being sectional semistable w.r.t.\
any of the numbers
$\overline{\delta}_1,...,\overline{\delta}_s$.
Since $(\E,\phi)$ is not semistable w.r.t.\ $\tilde{\delta}$, there
is either a saturated submodule $\E_0^\p\subset\ker(\phi)$ with
$\delta_{\E_0^\p}:=P_\E-(\rk\E/ \rk\E_0^\p)P_{\E_0^\p}<\tilde{\delta}$, or there
exists a saturated submodule $\E_0^{\p\p}\not\subset\ker(\phi)$ such that
$$\delta_{\E_0^{\p\p}}:=(\rk\E_0^{\p\p} P_\E- \rk\E
P_{\E_0^{\p\p}})/(\rk\E-\rk\E_0^{\p\p})>\tilde{\delta}\ .$$
In the first case suppose that
$\delta_{\E_0^\p}$ is minimal and in the second that
$\delta_{\E_0^{\p\p}}$ is maximal.
We consider only the first case, since the second can be treated similarly.
If $\delta_{\E_0^\p}\le 0$, then we are done.
Otherwise, set $\delta^\p_{i_0}:=\delta_{\E_0^\p}$.
By the minimality of $\delta_{i_0}^\p$, any submodule $\E^\p$ of
$\ker(\phi)$ satisfies
$$(\rk\E)\dim H^0(\E^\p(m)) \le \rk\E^\p (p-\overline{\delta^\p}_{i_0}),$$
and for $\E^\p=\E_0^\p$ we have equality.
Since $\E$ is not sectional semistable w.r.t.\ $\overline{\delta^\p}_{i_0}$,
there must exist a submodule $\E^{\p\p}\not\subset\ker(\phi)$ with
$$(\rk\E)\dim H^0(\E^{\p\p}(m))> \rk\E^{\p\p}
 (p-\overline{\delta^\p}_{i_0})+(\rk\E)\overline{\delta^\p}_{i_0}.$$
This makes it obvious that $(\E,\phi)$ cannot be sectional semistable w.r.t.\
to any positive rational number.\par
We still have to prove the ``stable'' version of the proposition.
For this we enlarge the constant $C$ such that $-C\le -\delta_i^0$, $i=1,...,s$,
where $\delta_i^0$ is the leading coefficient of $\delta_i$.
If $(\E,\phi)$ is a pair which is semistable w.r.t.\ the polynomial, say,
$\delta_{i_0}$ but not stable
w.r.t.\ any other polynomial $\delta$, then there must exist submodules
 $\E^\p\subset\ker(\phi)$
and $\E^{\p\p}$ belonging to $\tilde{\frak S}$ with
$${P_{\E^\p}\over\rk\E^\p}={P_{\E}-
\delta_{i_0}\over\rk\E}\q\hbox{and}\q {P_{\E^{\p\p}}-
\delta_{i_0}\over
\rk\E^{\p\p}}={P_{\E^\p}-\delta_{i_0}\over\rk\E}.$$
Since $m$ was so large that all modules in $\tilde{\frak S}$ are globally
 generated and
without higher cohomology, this gives
\begin{eqnarray*}
(\rk\E)\dim H^0(\E^\p(m)) &=& (\rk\E^\p)(p-\overline{\delta}_{i_0})\\
(\rk\E)\dim H^0(\E^{\p\p}(m)) &=& (\rk\E^{\p\p})(p-\overline{\delta}_{i_0})+
(\rk\E)\overline{\delta}_{i_0},\\
\end{eqnarray*}
and hence the assertion.
\end{pf}

\subsection{Proof of Theorem~\ref{MyStabCrit}}
\label{PfMyStabCrit}
For $b\in {\frak B}$, put $H_b:=\Hom(\bigwedge^rV,H^0(\L[\E_b](rm)))$ and
$\P_b:=\P((H_b\oplus S^rK)^\vee)$. The space
$\P_b$
admits the following natural $\C^*$-action:
$$z\cdot [h,\widehat{k}]:= [h, z\widehat{k}]=[z^{-1} h,\widehat{k}].$$
By \ref{C^*-ex}, this $\C^*$-action can be linearized
in such a way that the quotient is either $\P(H_b^\vee)$, $\P((S^rK)^\vee)$,
or
$\P(H_b^\vee)\times\P((S^rK)^\vee)$
equipped with the
polarization $[\O(a_1,a_2)]$ for any prescribed ratio $a_2/a_1$.
We are now able to apply our GIT-Theorem~\ref{GITThm}
to reduce Theorem~\ref{MyStabCrit}
to Theorem~\ref{StabCrit}.
\par
First we explain the assertion about semistability:
Suppose that $b=[q,\eps,\phi]$ lies in ${\frak B}^{ss}$.
Then either $\phi$ is injective, or $\E$ is semistable, or $\phi\neq 0$ and
the pair
$(\E,\phi)$ is semistable w.r.t.\ some $\delta_i$.
If $\phi$ is injective, we linearize in such a way that we obtain
$\P((S^rK)^\vee)$ as the quotient. By \ref{EasyStabCrit}, the point $[k]$ is
 semistable in $\P(K^\vee)$ and hence $[S^rk]$ is semistable in
$\P((S^rK)^\vee)$.
This implies by~\ref{GITThm} that $[h, S^rk]$ is semistable in $\P_b$.
If $\E$ is semistable, we linearize the $\C^*$-action in such a way that
the quotient $\P_b\catqot\C^*$ is given by $\P(H_b^\vee)$.
By \cite{Gi}, Theorem 0.7 (which does not depend on dimension 2),
the point $[h]$ is then semistable in $\P(H_b^\vee)$, and
hence $[h,S^rk]$ is $\SL(V)$-semistable in $\P_b$ by \ref{GITThm}.
If $\phi\neq 0$, $\eps\neq 0$ and $(\E,\phi)$ is semistable w.r.t.\
$\delta_i$, we choose the linearization of the $\C^*$-action in
such a way that the quotient is $\P(H_b^\vee)\times \P((S^rK)^\vee)$,
equipped with
a polarization $[\O(ra_1,a_2)]$ satisfying
$(a_2/a_1)= \rk\E \overline{\delta}_i/(p-\overline{\delta}_i)$.
By Theorem~\ref{StabCrit}, $([h], [S^rk])$ is semistable and thus $[h,S^rk]$
is semistable.
\par
Conversely, suppose $[h,S^rk]$ is $\SL(V)$-semistable.
By \ref{GITThm}
there is a linearization of the $\C^*$-action such that the image
of $[h,S^rk]$ is $\SL(V)$-semistable in the quotient
$\P_b\catqot\C^*$. There are three possible quotients:
If the quotient is $\P((S^rK)^\vee)$, then semistability implies that $[k]$
is semistable in $\P(K^\vee)$ and hence that $k$
is injective. It follows that $\E$ is a subsheaf of $\E_0$,
since we may assume that $m$ is so large that $\ker(\phi(m))$ is globally
generated.
If the quotient is $\P(H_b^\vee)$, then
 $\E$ is semistable by \cite{Gi}, loc.\ cit..
If the quotient is $\P(H_b^\vee)\times\P((S^rK)^\vee)$ with polarization
$[\O(a_1,a_2)]$, then $(\E,\phi)$ is sectional semistable w.r.t.\
$$\overline{\delta}:=p(ra_2/a_1)/(\rk\E +(ra_2/a_1))\ .$$

In view of \ref{SecStab} and \ref{TheProp}, $(\E,\phi)$ is semistable
w.r.t.\ some $\delta$,
hence $[q,\eps,\phi]$ lies in ${\frak B}^{ss}$.
\par
We still have to identify the stable points.
As the proof of \ref{GlueTog} shows, the oriented pair $(\E,\eps,\phi)$ given by
a point $b=[q,\eps,\phi]\in {\frak B}$ has only finitely many automorphisms if
and only if the associated point $[h,S^rk]\in\P_b$ has a finite
$\SL(V)$-stabilizer.
Let $b=[q,\eps,\phi]$ be a point whose associated element $[h,S^rk]$ in
$\P_b$ is stable.
If $h=0$ or $k=0$, then it is easy to see that the corresponding element
 $[S^rk]\in\P((S^rK)^\vee)$
or $[h]\in\P(H_b^\vee)$ is stable.
Hence $H^0(q(m))$ is an isomorphism and either $\phi$ is injective or $\E$
is a stable sheaf.
In both cases, the oriented pair $(\E,\eps,\phi)$ is stable and $H^0(q(m))$
 is an isomorphism,
in other words $b\in {\frak B}^s$.
If both $h\neq 0$ and $k\neq 0$, then by \ref{GITThm} $([h],[S^rk])\in
 \P(H_b^\vee)\times\P((S^rK)^\vee)$
is a polystable point w.r.t.\ the polarization, say, $\O(a_1,a_2)$.
By what we have already proved, $(\E,\phi)$ is a semistable pair.
Remark~\ref{properlysemistablepairs} shows that either $(\E,\phi)$ is a
stable pair
or there is an $i\in\{\, 1,...,s\,\}$ such that $(\E,\phi)$ is polystable
w.r.t.\ $\delta_i$.
In the first case, we are done.
In the second case, the finiteness of the stabilizer of $[h,S^rk]$ implies
that
the oriented pair $(\E,\eps,\phi)$ has only finitely many automorphisms, hence
it is a stable oriented pair.
\par
Suppose now that $b\in {\frak B}^s$.
If $\phi=0$, then $\E$ must be a stable coherent sheaf and thus
$[h]\in\P(H_b^\vee)$
is a stable point.
It follows that $[h,0]$ is a polystable point.
But as $[h,0]$ is a fixed point of the $\C^*$-action, the $\SL(V)$-stabilizer
 of $[h,0]\in\P_b$
can be identified
with the $\SL(V)$-stabilizer of $[h]\in \P(H_b^\vee)$,
so that $[h,0]$ is indeed
a stable point.
If $\eps=0$, then $\phi$ must be injective and we may argue in the same manner.
If both $\eps\neq 0$ and $\phi\neq 0$, it suffices to show that $[h,S^rk]$ is
a polystable point, since its stabilizer is finite by definition.
By the stability of $(\E,\phi)$, by the ``stable'' version of \ref{TheProp},
 and by the choice
of the $\delta_i$, there exists an index $i\in \{\, 1,...,s\,\}$
such that $(\E,\phi)$ is polystable w.r.t.\ $\delta_i$.
This in turn shows that $([h],[S^rk])\in \P(H_b^\vee)\times\P((S^rK)^\vee)$
 is polystable w.r.t.\
the linearization in $\O(ra_1,a_2)$ satisfying
$\overline{\delta}_i=p(a_2/a_1)/(\rk\E +(a_2/a_1))$.

\subsection{Moduli spaces of stable oriented pairs}
We need the following proposition
\begin{Prop}
\label{Proper}
The map ${\iota}_{\vert {\frak B}^{ss}}\colon {\frak B}^{ss}\lra \P({\frak
Z})^{ss}$  is
finite.
\end{Prop}
\begin{pf}
We claim that ${\iota}_{\vert {\frak B}^{ss}}$ is proper and injective.
Injectivity follows by standard arguments.
For the proof of properness, we will make use of the discrete valuative
criterion.
Let $C=\Spec R$ be the spectrum of a discrete valuation ring,
 $c_0\in C$ the closed point,
and $C_0:=C\setminus \{c_0\}$.
Suppose there is a commutative diagram:
\begin{center}
\unitlength=0.08mm
\begin{picture}(400,400)(0,60)
\put(0,390){$C_0$}
\put(75,400){\vector(1,0){220}}
\put(170,410){${\scriptstyle u}$}
\put(320,390){${\frak B}^{ss}$}
\put(20,350){\vector(0,-1){200}}
\put(340,350){\vector(0,-1){200}}
\put(0,90){$C$}
\put(70,100){\vector(1,0){220}}
\put(320,90){$\P({\frak Z})^{ss}$}
\put(170,120){$\scriptstyle \bar u$}
\put(380,240){$\scriptstyle \iota|_{{\frak B}^{ss}}$}
\multiput(70,140)(80,80){3}{\line(1,1){45}}
\put(250,320){\vector(1,1){36}}
\put(160,280){$\scriptstyle\tilde u$}
\end{picture}
\end{center}
We have to construct a lifting $\tilde{u}$ of the map $\overline{u}$.
By assumption, we are given a family
$({\frak E}_{C_0}, \eps_{C_0},\widehat{\phi}_{C_0},\O_{C_0})$
of semistable oriented pairs over
 $C_0\times X$. Note that ${\frak E}_{C_0}$ is
torsion free.
We claim that we can extend the quotient map
$$q_{C_0}\colon V\otimes \pi_X^*\O_X(-m)\lra
{\frak E}_{C_0}$$
to a homomorphism $q_C\colon V\otimes \pi_X^*\O_X(-m)\lra {\frak E}_C$
over $C\times X$,
where ${\frak E}_C$ is a flat family of torsion free coherent sheaves
 extending ${\frak E}_{C_0}$,
$q_{C}$ extends $q_{C_0}$, and its restriction to $\{c_0\}\times X$ is
generically
 surjective.
In order to prove this claim, we first extend the family ${\frak E}_{C_0}$
to a flat family
of quotients
$$\tilde{q}_C\colon V\otimes \pi_X^*\O_X(-m)\lra \tilde{\frak E}_C.$$
There is a locally free sheaf ${\cal H}$ on $X$ and an epimorphism
$\pi_X^*{\cal H}\lra \tilde{\frak E}^\vee_C$.
This yields a homomorphism
$$\lambda\colon V\otimes\pi_X^*\O_X(-m)\lra \tilde{\frak E}_C\lra \tilde {\frak
E}_C^{\vee\vee}\lra \pi_X^*{\cal H}^\vee.$$ Let ${\frak E}_C$ be the maximal
 subsheaf of
$\pi_X^*{\cal H}^\vee$ with the following properties
$${\frak E}_{C\vert C_0\times X}={\frak E}_{C_0};\qquad \Im\lambda\subset
{\frak
E}_C;\qquad
\dim(\mathop{\rm supp}({\frak E}_C/\Im\lambda))<\dim X.$$
Note that the set of subsheaves of $\pi_X^*{\cal H}^\vee$
having the above properties contains $\Im\lambda$.
One checks that ${\frak E}_{C\vert \{c_0\}\times X}$ is torsion free,
using arguments as in \cite{HL1}, p.85.
Let $t\in R$ be a generator of the maximal ideal.
There is a well defined integer $\alpha$ such that $(t^{\alpha}\eps_{C_0},
t^\alpha\widehat{\phi}_{C_0})$ extends to the family ${\frak E}_C$.

The classifying map to $\P({\frak Z})^{ss}$ induced by the resulting family
$$(q_C\colon V\otimes\pi_X^*\O_X(-m)\lra {\frak
E}_C,\tilde{\eps}_C,\tilde{\widehat{\phi}}_C,\O_C)$$
 is the same as the one induced by
$\overline{u}$. By the various stability criteria we have encountered so far,
 it follows
that $H^0(q_{C\vert \{c_0\}\times X}(m))$ is injective and that the triple
$({\frak E}_{C\vert \{c_0\}\times X}, \tilde{\eps}_{C\vert \{c_0\}\times X},
\tilde{\phi}_{c_0})$, where $\phi_{c_0}$ is a framing induced by
 $\tilde{\widehat{\phi}}_{C|\{c_0\}\times
X}$, is a semistable oriented pair.

Thus, ${\frak E}_{C\vert \{c_0\}\times X}(m)$ is globally generated
and without higher cohomology, the map $q_{C\vert \{c_0\}\times X}$ is
surjective,
and hence
$q_C\colon V\otimes \pi_X^*\O_X(-m)\lra {\frak E}_C$ is a
flat family of torsion free quotients.
The family $(q_C\colon V\otimes\pi_X^*\O_X(-m)\lra {\frak E}_C,\tilde{\eps}_C,
\tilde{\widehat{\phi}}_C,\O_C)$ defines by \ref{LocUnivProp} a morphism
$$\tilde{u}\colon C\lra {\frak B}^{ss}$$
which extends $u$ by construction.
\end{pf}

By Proposition 2.6.1. and \cite{Gi}, Lemma 4.6, the
quotient ${\frak B}^{ss}\catqot\SL(V)$  exists as a projective scheme. We  set
$${\cal M}_{(P,\L,\E_0)}^{ss}:={\frak B}^{ss}\catqot\SL(V)\ ,$$
$${\cal
M}_{(P,\L,\E_0)}^{s}:={\frak B}^s\catqot\SL(V)\ .$$

\begin{Thm}
\label{ModuliSpaces}
{\rm i)} There is a natural transformation of functors
$$ \tau\colon M_{(P,\L,\E_0)}^{ss}\lra h_{{\cal M}_{(P,\L,\E_0)}^{ss}},$$
\hspace*{0.5cm} such that for any scheme $\tilde{\cal M}$ and any natural
transformation
of functors
$$\tau^\p\colon M_{(P,\L,\E_0)}^{ss}\longrightarrow h_{\tilde{\cal M}}$$
\hspace*{0.5cm}there is a unique morphism $\vartheta\colon {\cal
M}_{(P,\L,\E_0)}^{ss}\lra
 \tilde{\cal M}$
such that
$\tau^\p=h(\vartheta)\circ \tau$.\\
\hspace*{0.5cm}{\rm ii)} The space ${\cal M}_{(P,\L,\E_0)}^s$ is a coarse moduli
 space for stable
oriented pairs.
\end{Thm}

\begin{pf}
The existence of the natural transformation is a direct consequence
of Proposition~\ref{LocUnivProp} and \ref{GlueTog}.
The minimality property of ${\cal M}_{(P,\L,\E_0)}^{ss}$  follows from
the universal property of the categorical quotient.\par
Since ${\frak B}^s$ is contained in the set of $\SL(V)$-stable points,
the set of closed points of ${\cal M}_{(P,\L,\E_0)}^s$ is the set of
equivalence classes of stable oriented pairs which means that
${\cal M}_{(P,\L,\E_0)}^s$ is a coarse moduli space.
\end{pf}
\vspace{0.3cm}

In our applications [OT2] we shall also need a slightly modified version of the
constructions and results above. We fix a line bundle ${\cal
L}_0\in\Pic (X)$ and consider  only torsion free sheaves of determinant
isomorphic to
${\cal L}_0$.

More precisely, an ${\cal L}_0$-{\it oriented pair of type
$(P,{\cal E}_0)$} is a triple $({\cal E},\varepsilon,\varphi)$ consisting of
a torsion free coherent sheaf ${\cal E}$ with Hilbert polynomial $P$ and
 with $\det{\cal
E}$ isomorphic to ${\cal L}_0$, a homomorphism
$\varepsilon:\det{\cal E}\map {\cal L}_0$,
and a homomorphism $\varphi:{\cal E}\map  {\cal E}_0$.

 Equivalence classes of such  ${\cal L}_0$-oriented pairs, families,
 equivalence classes
of families, (semi)stability and the corresponding functors
$M^{ss}_{(P,{\cal L}_0,{\cal
E}_0)}$ are defined as in 2.1. The same methods as above yield the
following result:

\begin{Thm}
\label{ModuliSpaces}
There exist moduli spaces ${\cal M}^{ss}_{(P,\L_0,\E_0)}$ and ${\cal
M}^{s}_{(P,\L_0,\E_0)}$ with the following properties:\\
\hspace*{0.5cm} {\rm i)} There is a natural transformation of
functors
$$ \tau\colon M_{(P,\L_0,\E_0)}^{ss}\lra h_{{\cal M}_{(P,\L_0,\E_0)}^{ss}},$$
\hspace*{0.5cm}  such that for any scheme $\tilde{\cal M}$ and any natural
transformation
of functors
$$\tau^\p\colon M_{(P,\L_0,\E_0)}^{ss}\longrightarrow h_{\tilde{\cal M}}$$
\hspace*{0.5cm} there is a unique morphism $\vartheta\colon {\cal
M}_{(P,\L_0,\E_0)}^{ss}\lra
\tilde{\cal
M}$ such that
$\tau^\p=h(\vartheta)\circ \tau$.\\
\hspace*{0.5cm} {\rm ii)} The space ${\cal M}_{(P,\L_0,\E_0)}^s$ is a
coarse moduli space
for
 stable
${\cal L}_0$-oriented pairs.
\end{Thm}

%%%%%%%%%%%%%%%%%%%%%%%%%%%%%%%%%%%%%%%%%%%%%%%%%%%%%%%%%%%%%%%%%%%%%%%%%%%%
%\section{Properties of the moduli spaces}

\subsection{The closed points of ${\cal M}^{ss}_{(P,\L,\E_0)}$}
Let $(\E,\eps,\phi)$ be a semistable oriented pair  of type $(P,\L,\E_0)$.
If $(\E,\eps,\phi)$ is stable, then it defines a closed point in  ${\cal
M}_{(P,\L,\E_0)}^{ss}$. If $(\E,\eps,\phi)$ is not stable, then either $\E$ is a
semistable but not stable coherent sheaf, or $\phi\neq 0$ and there exists a
$\delta\in\Q[x]$, $\delta>0$, such that
$(\E,\phi)$ is  semistable but not stable w.r.t.\ $\delta$.
In both cases, there is a  Harder-Narasimhan filtration
$$0= \E_0\subset \E_1\subset\cdots\subset\E_s=\E$$
of $\E$, whose associated graded sheaf $\mathop{\rm gr}(\E):=\bigoplus_{i=1}^s
\E_i/\E_{i-1}$ inherits a well-defined orientation $\eps_{\mathop{\rm gr}}$ and
a well-defined framing $\phi_{\mathop{\rm gr}}$
from $(\E,\eps,\phi)$.
As usual, the resulting object $({\mathop{\rm gr}}(\E),\eps_{\mathop{\rm gr}},
\phi_{\mathop{\rm gr}})$ is determined up to equivalence.
We call it the \it graded object associated to $(\E,\eps,\phi)$\rm .
Using the techniques of Section~\ref{PfMyStabCrit}, i.e.,  applying
\ref{GITThm} in the ``polystable'' version, we reduce the polystability of
$({\mathop{\rm gr}}(\E),\eps_{\mathop{\rm gr}},\phi_{\mathop{\rm gr}})$ to the
 respective
results of \cite{HL1}, \cite{HL2}, \cite{Gi}, and \cite{Ma2}.
Finally, one easily adapts the proof in \cite{HL2}, p.312, to show that a
semistable
oriented pair
$(\E,\eps,\phi)$ can be deformed into its graded object.
If we call two semistable oriented pairs $(\E_i,\eps_i,\phi_i)$, $i=1,2$,
\it gr-equivalent \rm if
their associated graded objects are equivalent, then
we see that the closed points of ${\cal M}^{ss}_{(P,\L,\E_0)}$
correspond to gr-equivalence classes of semistable oriented pairs of
type $(P,\L,\E_0)$.

\subsection{The $\C^*$-action on ${\cal M}^{ss}_{(P,\L,\E_0)}$}
\label{Flips}
The moduli space possesses a natural $\C^*$-action, given by
$$z\cdot [\E,\eps,\phi]:=[\E,\eps, z\phi]=[\E,z^{-r}\eps,\phi].$$
The set of fixed points of this action can easily be described:
It consists of classes $[\E,0,\phi]$, $[\E,\eps,0]$,
and of classes $[\ker(\phi)_{\max}\oplus\E^\p, \eps, \phi]$ with
$0\neq \ker(\phi)_{\max}$.
\par
The $\C^*$-action is naturally linearized in an ample line bundle
coming from the description of ${\cal M}^{ss}_{(P,\L,\E_0)}$
as GIT-quotient.
This line bundle and the polarization which it represents
may, however, depend on an integer $m$ chosen
in the course of the construction.
Nevertheless, we can state the following result which clarifies
the birational geometry of the moduli spaces ${\cal
M}^{ss}_{\delta}(X;\E_0,P)$ constructed in \cite{HL2}:
\begin{Thm}
Let $\delta_i\in\Q[x]$, $i=1,2$,
be polynomials with positive leading coefficients.
For a suitable choice of the polarization on ${\cal M}^{ss}_{(P,\L,\E_0)}$
the following properties hold true:
\par
{\rm i)}
${\cal M}^{ss}_{\delta_i}(X;\E_0,P)$, $i=1,2$, are $\C^*$-quotients of
the master space ${\cal M}^{ss}_{(P,\L,\E_0)}$.
\par
{\rm ii)}
${\cal M}^{ss}_{\delta_1}(X;\E_0,P)$ and ${\cal M}^{ss}_{\delta_2}(X;\E_0,P)$
are related by a chain of generalized flips.
\end{Thm}
\begin{pf}
Let $m$ be so large that a pair $(\E,\phi)$ is
semistable w.r.t.\ $\delta_i$ if and only the pair $(\E(m),\phi(m))$
is sectional semistable w.r.t.\ $\delta_i(m)$, $i=1,2$, and that
all the other requirements needed in the constructions are met.
Then our proof of Theorem~2.4.1 together with the results
of Section~1 easily yields the assertions of the theorem.
\end{pf}
We note that the $\delta_i$ for which the corresponding
set of $\C^*$-stable points meets the fixed point set of the
$\C^*$-action, i.e.,
for which the corresponding set of $\C^*$-stable points
contains stable oriented pairs of the type
$[\ker(\phi)_{\max}\oplus\E^\p, \eps, \phi]$ with
$0\neq \ker(\phi)_{\max}$ are uniquely determined.
The corresponding polynomial is
$\rk\E^\p(P_{\E^\p}-P_{\ker(\phi)_{\max}}/\rk\ker(\phi)_{\max})$.
The associated moduli spaces ${\cal M}_{\delta_i}$ are
those which show up ``at the top'' of the flips.
\newpage

\vspace{1cm}
\noindent
e-mail: {\small okonek@@math.unizh.ch, schmitt@@math.unizh.ch, teleman@@math.unizh.ch}

\end{document}